%% file: 2003-1.tex
\documentclass{gtart}

\input gtoutput

\lognumber{253}
\volumenumber{7}\papernumber{1}\volumeyear{2003}
\pagenumbers{1}{31}
\received{9 May 2002}
\accepted{8 November 2002}
\published{23 January 2003}
\proposed{Vaughan Jones}
\seconded{Yasha Eliashberg, Joan Birman}

\usepackage{amssymb,amsmath,dptsymb,graphics,subfigure,calc}

\newcommand{\RR}{\mathbb R}

\newcommand{\ZZ}{\mathbb Z}
\newcommand{\comma}{\mathbin ,}
\newcommand{\sltwo}{{{\mathfrak sl}_2}}
\newcommand{\frakg}{{\mathfrak g}}
\newcommand{\isom}{\simeq}
\newcommand{\nisom}{\not\simeq}

\newcommand{\bigcircle}{\bigcirc}

\newcommand{\tensor}{\otimes}

\newcommand{\ad}{\operatorname{ad}}

\newcommand{\tr}{\operatorname{tr}}
\renewcommand{\det}{\operatorname{det}\nolimits}

\newcommand\strutn[2]{{\,^{#1}\!\!\frown^{#2}}}
\newcommand\strutu[2]{{\,_{#1}\!\!\smile_{#2}}}
\newcommand\strutv[2]{{\mathop\mid^{#1}_{#2}}}
\newcommand{\Strutn}{{\!\frown}}

\newcommand{\connect}{\mathbin \#}
\newcommand{\cA}{{\mathcal A}}
\newcommand{\cB}{{\mathcal B}}

\newcommand{\directedcircle}{{\circlearrowright}}

\theoremstyle{plain}
\newtheorem{theorem}{Theorem}
\newtheorem{proposition}{Proposition}[section]
\newtheorem{lemma}[proposition]{Lemma}
\newtheorem{corollary}[proposition]{Corollary}

\newtheorem*{WHEELING}{Wheeling Theorem} 
\newtheorem*{WHEELS}{Wheels Theorem} 

\theoremstyle{definition}

\theoremstyle{remark}

\newtheorem{exercise}[proposition]{Exercise}
\newtheorem{hint}[proposition]{Hint}
\newtheorem{remark}[proposition]{Remark}

\DeclareMathOperator{\Vac}{(Vacuum)}

\DeclareMathOperator{\CCable}{C\Delta}

\newcommand{\Zero}{0}
\newcommand{\One}{1}
\newcommand{\Two}{2}
\newcommand{\hopf}[2]{{\,_{#1}\!\Hopf_{#2}}}
\newcommand{\Hopf}{\HopfLink}
\newcommand{\openhopf}[2]{\OpenHopf_{\!#1}^{#2}}

\newcommand{\bc}{{\text{bc}}}

\def\lbl#1{\label{#1}}

\newcommand{\fig}[2]{%
\input{draws/#1.tex}%
}
\newcommand{\figcent}[2]{\mathcenter{%
\input{draws/#1.tex}%
}}

\def\mathcenter#1{
	\raisebox{.5ex+\depth-.5\totalheight}{\hbox{#1}}
}

\newcommand{\mlarge}[1]{{\text{\Large $#1$}}}

\def\:{\kern .1ex\colon\thinspace}

\begin{document}

\title[Applications of knot theory to Lie algebras and Vassiliev invariants]{
  Two applications of elementary knot theory to\\Lie algebras and
  Vassiliev invariants
}

\authors{Dror Bar-Natan\qquad\quad  Thang T\thinspace Q Le\qquad\quad Dylan P Thurston}
\coverauthors{Dror Bar-Natan\\Thang T\noexpand\thinspace Q Le\\Dylan P Thurston}
\asciiauthors{Dror Bar-Natan\\Thang T Q Le\\Dylan P Thurston}
\shortauthors{Bar-Natan, Le and Thurston}
\address{
  \parbox{1.7in}{\centering
    Dept of Mathematics\\
    University of Toronto\\
    Toronto ON M5S 3G3\\
    Canada
  }
  \parbox{1.7in}{\centering
    Dept of Mathematics\\
    SUNY at Buffalo\\
    Buffalo NY 14214\\
    USA
  }
  \parbox{1.7in}{\centering
    Dept of Mathematics\\
    Harvard University\\
    Cambridge, MA 02138\\
    USA
  }\\\bigskip\\
{\tt drorbn@math.toronto.edu\quad letu@math.buffalo.edu\quad dpt@math.harvard.edu}
\\\smallskip\\
{\tt  http://www.math.toronto.edu/\char'176drorbn\qua 
buffalo.edu/\char'176letu\qua harvard.edu/\char'176dpt}
}
\asciiaddress{Department of Mathematics, University of Toronto\\Toronto 
ON M5S 3G3, Canada\\ Department of Mathematics, SUNY at Buffalo\\Buffalo 
NY 14214, USA\\Department of Mathematics, Harvard University\\Cambridge, 
MA 02138, USA}
\asciiemail{drorbn@math.toronto.edu, letu@math.buffalo.edu,
dpt@math.harvard.edu}

\primaryclass{57M27}

\secondaryclass{17B20, 17B37}

\keywords{Wheels, Wheeling, Vassiliev invariants, Hopf link, $1+1=2$, Duflo
isomorphism, cabling}

\begin{abstract}
Using elementary equalities between various cables of the unknot
and the Hopf link, we prove the Wheels and Wheeling conjectures
of~\cite{BGRT:WheelsWheeling, Deligne:letter}, which give,
respectively, the exact Kontsevich integral of the unknot and a
map intertwining two natural products on a space of diagrams. It
turns out that the Wheeling map is given by the Kontsevich
integral of a cut Hopf link (a bead on a wire), and its
intertwining property is analogous to the computation of $1+1=2$
on an abacus. The Wheels conjecture is proved from the fact that
the $k$-fold connected cover of the unknot is the unknot for all
$k$.

Along the way, we find a formula for the invariant of the general
$(k,l)$ cable of a knot. Our results can also be interpreted as a
new proof of the multiplicativity of the Duflo--Kirillov map
$S(\frakg) \to U(\frakg)$ for metrized Lie
(super-)algebras $\frakg$.
\end{abstract}
\asciiabstract{Using elementary equalities between various cables of
the unknot and the Hopf link, we prove the Wheels and Wheeling
conjectures of [Bar-Natan, Garoufalidis, Rozansky and Thurston,
arXiv:q-alg/9703025] and [Deligne, letter to Bar-Natan, January 1996,
http://www.ma.huji.ac.il/\char'176drorbn/Deligne/], which give,
respectively, the exact Kontsevich integral of the unknot and a map
intertwining two natural products on a space of diagrams. It turns out
that the Wheeling map is given by the Kontsevich integral of a cut
Hopf link (a bead on a wire), and its intertwining property is
analogous to the computation of 1+1=2 on an abacus. The Wheels
conjecture is proved from the fact that the k-fold connected cover of
the unknot is the unknot for all k. Along the way, we find a formula
for the invariant of the general (k,l) cable of a knot. Our results
can also be interpreted as a new proof of the multiplicativity of the
Duflo-Kirillov map S(g)-->U(g) for metrized Lie (super-)algebras g.}

\maketitlepage

\section{Introduction}

\subsection{The Duflo--Kirillov isomorphism} \lbl{sec:intro-Duflo}   The
Duflo--Kirillov isomorphism is an {\em algebra} isomorphism between
the invariant part of the symmetric algebra and the center of the
universal enveloping algebra for any Lie algebra $\frakg$. This
isomorphism was first described for semi-simple Lie algebras by
Harish-Chandra. Kirillov gave a formulation of the Harish-Chandra
map that has meaning for all finite-dimensional Lie algebras, and
conjectured that it is always an algebra isomorphism. The
conjecture was proved by Duflo~\cite{Duflo:Operateurs}. Although
the Kirillov--Duflo map can be formulated in a very explicit way as
a linear map between two pretty simple algebras (with very
explicit structure), all known proofs of the Duflo theorem were
difficult: In the book of Dixmier \cite{Dixmier}, the proof is
given only in the last chapter and it utilizes most of results
developed in the whole book, including many classification
results (a situation Godement \cite{Godement} called
``scandalous"). As discussed below, there have been several
recent proofs that do not use classification results, but they
all use tools from well outside the natural domain of the problem.

Let us review briefly the Duflo theorem. The
Poincar{\'e}--Birkhoff--Witt map between the symmetric algebra and the
universal enveloping algebra of a Lie algebra $\frakg$,
\[
\chi\: S(\frakg) \longrightarrow U(\frakg),
\]
given by taking a monomial $x_1\dots x_n$ in $S(\frakg)$ and
averaging over the product (in $U(\frakg)$) of the $x_i$ in all
possible orders, is an isomorphism of vector spaces and
$\frakg$-modules. Since $S(\frakg)$ is abelian and $U(\frakg)$ is
generally not, $\chi$ is clearly not an algebra isomorphism. Even
restricting to the invariant subspaces on both sides,
\[
\chi\: S(\frakg)^\frakg \longrightarrow U(\frakg)^\frakg =
\text{center of $U(\frakg)$},
\]
$\chi$ is still not an isomorphism of algebras.

The Duflo theorem says that the combination $\chi \circ
\partial_{j^{\frac12}}$, with $\partial_{j^{\frac12}}\: S(\frakg)
\longrightarrow S(\frakg)$  defined below, is an algebra
isomorphism between $ S(\frakg)^\frakg$ and $U(\frakg)^\frakg$.

Here $j^{\frac12}(x)$ is a formal power series (beginning with 1)
on $\frakg$ , defined by
$$j^{\frac12}(x) =
 \det\nolimits^{\frac12}\left(
  \frac{\sinh(\frac{1}{2} \ad x)}
       {\frac{1}{2} \ad x}\right).$$
The operator $\partial_{j^{\frac12}}$ is obtained by plugging the
(commuting) vector fields $\partial/\partial x^*$ (on $\frakg^*$)
in the power series $j^{\frac12}$. (Note that for
$x^*\in\frakg^*$, $\partial/\partial x^*$ transforms like an
element of $\frakg$).  The result is an infinite-order
differential operator on $\frakg^*$, which we can then apply to a
polynomial on $\frakg^*$ ($\equiv$ an element of $S(\frakg)$).
For details, see \cite{Duflo:Operateurs}. The function
$j^\frac{1}{2}(x)$ plays an important role in Lie theory.
  Its square, $j(x)$, is the Jacobian of the
exponential mapping from $\frakg$ to the Lie group $G$. The
operator $\partial_{j^{\frac12}}$ is called {\em the strange
isomorphism} by Kontsevich
\cite{Kontsevich:DeformationQuantization}.

\subsection{Elementary knot theory}\lbl{sec:ElemKnotTheory} We will
touch upon two  simple facts in knot theory that have deep
consequences for Lie algebras and Vassiliev invariants.  The two
facts can be summarized by the catch phrases ``$1+1=2$'' and
``$n\cdot 0 = 0$.''

\begin{itemize}
\item ``$1+1=2$.'' This refers to a fact in ``abacus
arithmetic.''  On an abacus, the number $1$ is naturally
represented by a single bead on a wire, as in
Figure~\ref{fig:One}, which we think of as a tangle.  The fact
that $1+1 = 2$ then becomes the equality of the two tangles in
Figure~\ref{fig:OneOneTwo}.  On the left side of the figure,
``$1+1$'', the two beads are well separated, as for connect sum
of links or multiplication of tangles; on the right side, ``$2$'',
we instead start with a single bead and double it, so the two
beads are very close together.

In other terms, the connected sum of two Hopf links is the same as
doubling one component of a single Hopf link, as in
Figure~\ref{fig:HopfSum}.
\begin{figure}[ht!]%
\begin{center}%
  \subfigure[The link ``$\One$'']%
   {\makebox[.25\linewidth]{%
\input{draws/One.tex}%
\lbl{fig:One}}}%
\hfil\subfigure[``$\One + \One = \Two$'']%
   {\makebox[.35\linewidth]{%
\input{draws/OOneTwo.tex}%
\lbl{fig:OneOneTwo}}}%
\hfil\subfigure[An alternate version of ``$\One + \One = \Two$'']%
   {\makebox[.35\linewidth]{%
\input{draws/HopfSum.tex}%
\lbl{fig:HopfSum}}}%
\end{center}%
\caption{Elementary knot theory, part 1}\lbl{fig:Abacus}%
\end{figure}
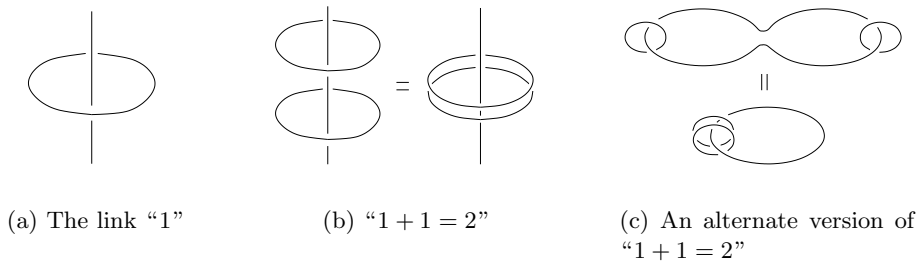

\item ``$n\cdot 0 = 0.$''  In the spirit of abacus
arithmetic, $0$ is represented as just a single vertical strand.
We prefer to close it off, yielding the knot in
Figure~\ref{fig:Zero}. The knot $n\cdot 0$ is then this knot
repeated $n$ times, as in
Figure~\ref{fig:nZero}.%
The two knots are clearly the same, up to framing.

\begin{figure}[ht!]%
\begin{center}%
  \subfigure[The knot ``$\Zero$'']%
   {\makebox[.3\linewidth]{%
\input{draws/Zero.tex}%
\lbl{fig:Zero}}}%
\qquad\subfigure[The knot ``$n\cdot\Zero$,'' shown here for $n=3$]%
   {\makebox[.3\linewidth]{%
\input{draws/nZero.tex}%
\lbl{fig:nZero}}}%
\end{center}%
\caption{Elementary knot theory, part 2}%
\end{figure}
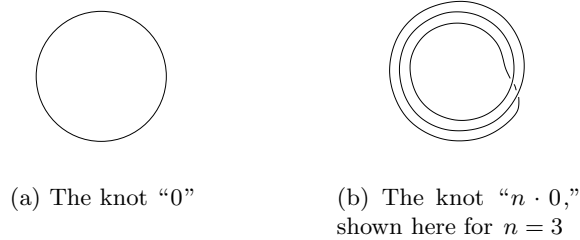
\end{itemize}

\subsection{Wheels and wheeling: main results} The bridges between the knot theory of
Section~\ref{sec:ElemKnotTheory} and the seemingly quite
disparate Lie algebra theory of Section~\ref{sec:intro-Duflo} are
a certain spaces of uni-trivalent diagrams (called Jacobi
diagrams) modulo local relations.  (See Section
\ref{Preliminaries}, the 1-valent vertices are called the
``legs'' of the diagram.)  On the one hand, such diagrams give
elements of $U(\frakg)$ or $S(\frakg)$ for every metrized Lie
algebra $\frakg$ in a uniform way; on the other hand, they occur
naturally in the study of finite type invariants of knots
\cite{BarNatan:OnVassiliev,Kontsevich}. Like the associative
algebras $S(\frakg)$ and $U(\frakg)$ associated to Lie algebras,
these diagrams appear in two different varieties: $\cA$, in which
the legs have a linear order, as in Figure~\ref{fig:A-examp}, and
$\cB$, in which the legs are unordered, as in
Figure~\ref{fig:B-examp}.  As for Lie algebras, they each have a
natural algebra structure (concatenation and disjoint union,
respectively); and there is an isomorphism $\chi\:\cB\to
\cA$ between the two (averaging over all possible orders of the
legs).

\begin{figure}[ht!]%
\begin{center}%
  \subfigure[A sample element of $\cA$]%
   {\makebox[.35\linewidth]{%
\input{draws/A-examp.tex}%
\lbl{fig:A-examp}}}%
\qquad\subfigure[A sample element of $\cB$]%
   {\makebox[.35\linewidth]{%
\input{draws/B-examp.tex}%
\lbl{fig:B-examp}}}%
\end{center}%
\caption{Examples of Jacobi diagrams}%
\end{figure}
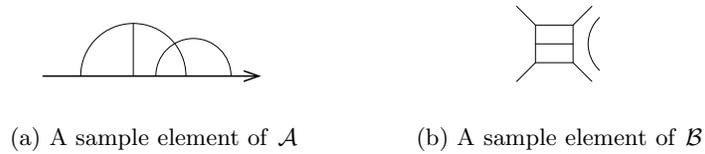

There is one element of the algebra $\cB$ that will be
particularly important for us: the ``wheels'' element.  It is the
diagrammatic analogue of the function $j^\frac12$ above:
\begin{equation} \lbl{eq:OmegaDef}
  \Omega=\exp \sum_{n=1}^\infty b_{2n}\omega_{2n} \in \cB,
\end{equation}
where:
\begin{itemize}
\item The `modified Bernoulli numbers' $b_{2n}$ are defined by the power
series expansion
\begin{equation} \lbl{eq:MBNDefinition}
  \sum_{n=0}^\infty b_{2n}x^{2n} = \frac{1}{2}\log\frac{\sinh x/2}{x/2}.
\end{equation}
These numbers are related to the usual Bernoulli numbers $B_{2n} =
4n\cdot(2n)!\cdot b_{2n}$ and to the values of the Riemann
$\zeta$-function on the even integers. The first three modified
Bernoulli numbers are $b_2=1/48$, $b_4=-1/5760$, and
$b_6=1/362880$.

\item The `$2n$-wheel' $\omega_{2n}$ is the degree $2n$ Jacobi diagram
made of a $2n$-gon with $2n$ legs:
\begin{equation}\label{eq:omega-definition}
  \omega_2=\mathcenter{%
\input{draws/2wheel.tex}},\quad
  \omega_4=\mathcenter{%
\input{draws/4wheel.tex}},\quad
  \omega_6=\mathcenter{%
\input{draws/6wheel.tex}},\quad\ldots.
\end{equation}
\end{itemize}

Let $\partial_\Omega$ be the operation of applying $\Omega$ as a
differential operator, which takes a diagram $D$ and attaches
some of its legs to all the legs of $\Omega$.  (See
Section~\ref{sec:diagram-differential} for the precise
definition.)

The first main result of this paper is the following analog of the
Duflo theorem.
\begin{WHEELING}
The map $\Upsilon = \chi \circ \partial_\Omega\: \cB \to
\cA$ is an algebra isomorphism.  \end{WHEELING}

Although the Wheeling theorem was motivated by Lie algebra
considerations when it was first
conjectured~\cite{BGRT:WheelsWheeling, Deligne:letter}, the proof
we will give, based on the equation ``$1+1=2$'' from
Section~\ref{sec:ElemKnotTheory}, is entirely independent of Lie
algebras and is natural from the point of view of knot theory. In
particular, we obtain a new proof of the Duflo theorem for
metrized Lie algebras, with some advantages over the original
proofs by Harish-Chandra, Duflo, and Cartan: our proof does not
require any detailed analysis of Lie algebras, and so works in
other contexts in which there is a Jacobi relation.  For
instance, our proof works for super Lie algebras with no
modification.

The Wheeling theorem has already seen several applications.  We
will use it to compute the Kontsevich integral of the unknot,
using our second elementary knot theory identity ``$n\cdot 0 =
0$''.

\begin{WHEELS}
The Kontsevich integral of the unknot is
\[
Z(\bigcircle) = \Omega \in \cB.
\]
\end{WHEELS}

The Wheeling  theorem was first conjectured by
Deligne~\cite{Deligne:letter} and by Bar-Natan, Garoufalidis,
Rozansky and Thurston~\cite{BGRT:WheelsWheeling}, who also
conjectured the Wheels theorem.

Along the way we also find a formula describing the behaviour of
the Kontsevich integral under connected cabling of knots. We also
compute the Kontsevich integral of the Hopf link $\Hopf$; which is
intimately related to the map $\Upsilon$ above.

Further computations for a sizeable class of knots, links, and
3-manifolds (including all torus knots and Seifert-fiber homology
spheres) have been done by Bar-Natan and
Lawrence~\cite{Bar-NatanLawrence:RationalSurgery}.  Hitchin and
Sawon~\cite{HitchinSawon:CurvatureHyperKahler} have used the
Wheeling theorem to prove an identity expressing the $L^2$ norm of
the curvature tensor of a hyperk{\"a}hler manifold in terms of
Pontryagin classes.  In a future paper~\cite{Thurston:Torus} one
of us (DPT) will show how to write simple formulas for the
action of $\sltwo(\ZZ)$ on the vector space associated to a torus
in the perturbative TQFT of Murakami and
Ohtsuki~\cite{MO:TQFTUniversal}. Our connected cabling formula
also finds application in recent work of Roberts and Willerton on
the ``total Chern class" invariant of knots.

There are two other recent proofs of the Wheeling theorem. One is
due to
Kontsevich~\cite[Section~8]{Kontsevich:DeformationQuantization},
as expanded by~\cite{ADS:KontsevichQuantization,HV,MT}.
Kontsevich's proof is already at a diagrammatic level, similar to
the one in this paper, although it is more general: it works for
all Lie algebras, not just metrized ones.  His proof again uses a
transcendental integral,
  similar in spirit to the ``Kontsevich integral" in the theory of
  Vassiliev invariants \cite{Kontsevich}. Another
proof is due to Alekseev and
Meinrenken~\cite{AlekseevMeinrenken:NonCommutativeWeil}.  The
Alekseev and Meinrenken paper is not written in diagrammatic
language, but seems to extend to the diagrammatic context without
problems.  Their proof does not involve transcendental integrals:
the only integral in their proof is in the proof of the Poincar{\'e}
lemma (the homology of $\RR^n$ is trivial in dimension $>\,0$).

\subsection{Plan of the paper} In the first section  we
review the theory of Jacobi diagrams.  Section \ref{cabling} is
devoted to  cabling formulas of the Kontsevich integral which are
crucial in the proofs of main theorems. In Sections~\ref{wheeling}
and \ref{wheels} we prove the Wheeling and Wheels theorems. In
Section \ref{Hopf} we calculate the values of the Kontsevich
integral of the Hopf link. In the Appendix we give a
self-contained method to determine the coefficients of the wheels
element.

\subsection{Acknowledgement} Research by the authors DBN and DPT was
supported in part by BSF grant \#1998-119. The author TTQL was partially
supported by NSF grant DMS-9626404 and a postdoc fellowship at the
Mathematical Sciences Research Institute in Berkeley in 1996--1997.
Research at MSRI was supported in part by NSF grant DMS-9022140. The
author DPT was supported by an NSF Graduate Student Fellowship, a
Sloan Dissertation Fellowship, and an NSF Postdoctoral Research
Fellowship. In addition we wish to thank A~Haviv, J~Lieberum,
A~Referee and J~Roberts for their comments and suggestions.

\section{Preliminaries on Jacobi diagrams} We recall basic definitions
and some known properties of Jacobi diagrams in this section. For
details, see \cite{BarNatan:OnVassiliev}.

\subsection{Jacobi diagrams}\label{Preliminaries}
 An {\em open Jacobi diagram} (sometimes called a
Chinese Character, uni-trivalent graph, or web diagram) is a
vertex-oriented uni-trivalent graph, ie, a graph with univalent
and trivalent vertices together with a cyclic ordering of the
edges incident to the trivalent vertices.  Self-loops  and
multiple edges are  allowed.  A univalent vertex is called {\em a
leg}, and trivalent vertex is also called an {\em internal
vertex}. In planar pictures, the orientation on the edges
incident on a vertex is the clockwise orientation, unless
otherwise stated. The {\em degree} of an open Jacobi diagram is
half the number of vertices (trivalent and univalent). Some
examples are shown in Figure~\ref{fig:B-examp}.

Suppose  $X$ is a compact oriented 1-manifold (possibly with boundary,
often with labeled components) and $Y$ a finite set of (labeled)
asterisks, symbols of the form $\ast_x$, $\ast_y$, etc.. A {\em Jacobi
diagram based on $X\cup Y$} is a graph $D$ together with a
decomposition $D = X \cup \Gamma$, where $\Gamma$ is an open Jacobi
diagram with some legs labeled by elements of $Y$, such that $D$ is the
result of gluing all the non-labeled legs of $\Gamma$ to distinct
interior points of $X$. Note that repetition of labels is allowed, and
not all labels have to be used. The {\em degree} of $D$, by definition,
is the degree of $\Gamma$. Usually $X$ is called the {\em skeleton} of
$D$, and in picture is depicted by bold lines.

 Suppose
$\phi\: X'\to X$ is a covering map between compact oriented
1-manifolds, and $D= X\cup \Gamma$ is a Jacobi diagram based on
$X\cup Y$. The {\em pull-back} $\phi^*(D)$ is the sum over all
Jacobi diagrams $D'$ based on $X'\cup Y$ such that $\phi(D')=D$.
Here $\phi(D')= D$ means $D'=X'\cup \Gamma$ and $\phi$ can be
extended to $D'$ so that it is identity on $\Gamma$.

The space $\cA^f(X\cup Y)$, $X$ and $Y$ as above, is the space of
Jacobi diagrams based on $X\cup Y$ modulo the usual antisymmetry,
IHX and STU relations (see \cite{BarNatan:OnVassiliev}). The
completion of $\cA^f(X\cup Y)$ with respect to degree is denoted
by $\cA(X\cup Y)$. When $\phi\:X' \to X$ is a cover, the pull-back
$\phi^*$ descends to a well-defined map from $\cA(X\cup Y)$ to
$\cA(X'\cup Y)$. An example of pull-backs is the Adams operation
in \cite{BarNatan:OnVassiliev}.

Let $\cA^\bc(X\cup Y)$ be the subspace of $\cA(X\cup Y)$ spanned
by boundary-connected Jacobi diagrams: diagrams with no connected
components that are disjoint from the skeleton $X$.

There is a natural map from $\cA(\uparrow\cup X)$ to
$\cA(\directedcircle \cup X)$ given by attaching the two
endpoints of the interval $\uparrow$.  If $X$ is a closed
1-manifold, then  this map is an isomorphism. In particular, when
$X=\emptyset$, the spaces $\cA(\uparrow)$ and
$\cA(\directedcircle)$ are canonically isomorphic. But this is
not true  if $X$ has an interval component. Explicitly,
$\cA(\uparrow\uparrow)\nisom\cA(\uparrow\directedcircle)\isom\cA(\directedcircle\directedcircle)$.

An open  Jacobi diagram  is {\em strutless} if it does not have a
connected component homeomorphic to a strut $\strutn{}{}$, ie\ an
interval. A {\em strutless element} of $\cA(Y)$, where $Y$ is a
set of asterisks, is a linear combination of strutless diagrams.

\subsection{Special interesting cases} Of special interest are the
following $\cA(X\cup Y)$.

For $X=\emptyset$ and $Y$ has one element, the space $\cA(X\cup Y)$ is
denoted by $\cB$. Note that all the labels of legs of diagrams in
$\cB$ are the same, and we often forget the labels. There is a
natural product in $\cB$ defined by taking disjoint union of
diagrams. With this product $\cB$ is a commutative algebra. The
wheels $\omega_{2n}$ introduced in the introduction belongs to
$\cB$.

For $X=\directedcircle$, the oriented circle, and $Y=\emptyset$,
the space $\cA (\directedcircle)$, also denoted simply by $\cA$,
is the space in which lie the values of the Kontsevich integral
of a knot. There is a natural product in $\cA$ defined by taking
connected sums of diagrams based on $\directedcircle$. With this product
$\cA$ is a commutative algebra. As noted before, $\cA$ is
canonically isomorphic to $\cA(\uparrow)$, and we will often
identify these vector spaces. Note that the space $\cA$ of
\cite{BarNatan:OnVassiliev, LeMurakami:Parallel} is equal to our
$\cA^\bc(\directedcircle)$, the boundary-connected part.

Suppose $X=Y=\emptyset$. The space $\cA(\emptyset)$ is the space
in which lie the values of the LMO invariants of 3-manifolds
\cite{LMO}. With disjoint union as the product, $\cA(\emptyset)$
becomes a commutative algebra, and all other $\cA(X\cup Y)$ have a
natural $\cA(\emptyset)$-module structure.

It is known that for any metrized Lie algebra $\frakg$, there are the
weight maps, which are algebra homomorphisms, $W_{\frakg}\:\cB^f \to
S(\frakg)^\frakg$ and $W_{\frakg}\:  \cA^f(\directedcircle) \to U(\frakg)^\frakg$,
see \cite{BarNatan:OnVassiliev}.  Here $S(\frakg)$ and $U(\frakg)$ are
respectively  the symmetric algebra and the universal enveloping
algebra of $\frakg$, and $M^\frakg$ is the invariant subspace of the
$\frakg$-module $M$.  Thus $U(\frakg)^\frakg$ is the center of
$U(\frakg)$. In some sense, one can think of $\cA$ and $ \cB$ as being
related to a ``universal (metrized) Lie algebra'', incorporating
information about all Lie algebras at once. But $\cB$ and $\cA$ are
both bigger and smaller than that. For example, the map from $\cB$ to
the product of $S(\frakg)^\frakg$ for all metrized Lie algebras is
neither injective nor surjective: There are elements of $\cB$ that are
non-zero but become zero when evaluated in any metrized Lie algebra
\cite{Vogel:Structures, Lieberum:NotComing}\footnote{These references
only deal with semi-simple Lie (super-) algebras, but according to
Vogel and Lieberum (via private communications), Vogel's results extend
to all metrized Lie (super-) algebras.}.  Not all elements of
$S(\frakg)^\frakg$ are in the image of the map $W_\frakg$. (For
instance, the image of $W_\frakg$ consists of polynomials of even order
only.)

\subsection{Symmetrization maps}
 One can define an analog of the Poincare--Birkhoff--Witt
 isomorphism
for diagrams as follows.

Suppose $X$ is a collection of compact oriented 1-manifolds and
asterisks. The symmetrization map $\chi_x\:\cA(\ast_x\cup X)
\to\cA(\uparrow_x\cup X)$ is a linear map  defined on a
diagram $D$ by taking the average over all possible ways of
ordering the legs labeled by $x$ and attach them to an oriented
interval. It is known that $\chi_x$ is a vector space isomorphism
\cite{BarNatan:OnVassiliev}.

In particular, the symmetrization map $\chi\: \cB \to
\cA(\uparrow) \equiv \cA(\directedcircle)$ is an isomorphism of vector
spaces, but it is not an algebra isomorphism. We drop the label
here. The two products, disjoint union and connected sum, live on
isomorphic spaces $\cB$ and $\cA$, and may be confused. We
usually write out the product in cases of ambiguity.

\subsection{Symmetrization for closed components of the skeleton}
\label{sec:diagrams2} We have seen that, using the symmetrization
map, one can trade an oriented interval in $X$ with an asterisk.
We want to do the same with closed component in $X$. For this we
need the {\em link relations}.

Suppose  $\ast_y$ is an element of $Y$. If a leg of a diagram is
labeled $y$, then the edge having this leg as an end is called
{\em a $y$-edge}. In $\cA(X\cup Y)$, {\em link relations} on $y$
are parametrized by Jacobi diagrams based on $X\cup Y$ in which
one of the $y$-labeled legs is distinguished. The corresponding
link relation is the sum of all ways of attaching the
distinguished leg to all the other $y$-edges:
\[
\figcent{l-rel-0}{0.4} \mapsto \figcent{l-rel-1}{0.4} =
  \figcent{l-rel-2}{0.4}
+ \figcent{l-rel-3}{0.4} + \cdots + \figcent{l-rel-4}{0.4}.
\]
Suppose $X$ is a compact oriented 1-manifold, $Y$ is a set of asterisks
$\ast$, and $Y'$ is a set of circled asterisks, symbols of the form
$\circledast_x$, $\circledast_y$, etc. Define
 $\cA(X \cup Y \cup Y')$ as the space of Jacobi diagrams based on $X
\cup Y \cup Y'$ modulo  the anti-symmetry, IHX, and STU relations as
before and, in addition, link relations on each label in $Y'$.

Suppose a circled asterisk $\circledast_y$ is not in $Y'$. The
symmetrization map $\chi_x\:\cA(\circledast_y\cup X\cup Y\cup Y')
\to\cA(\directedcircle_y\cup X\cup Y\cup Y')$ is the
linear map defined on a diagram $D$ by taking the average over all
possible ways of cyclic-ordering the legs labeled by $y$ and
attach them to the circle $\directedcircle_y$. It is known that
$\chi_y$ is a vector space isomorphism \cite{BGRT2}.

\subsection{Diagrammatic Differential Operators}
\label{sec:diagram-differential}
 For a strutless diagram $C
\in \cB$, the operation of {\em applying C as a differential
operator}, denoted $\partial_C\: \cB \to \cB$, is defined
to be
\[ \partial_C(D)=\begin{cases}
    0 & \parbox{1.7in}{if $C$ has more legs than $D$,} \\[3pt]
    \parbox{2.2in}{
      the sum of all ways of gluing all the legs of $C$ to some (or all)
      legs of $D$
    }\quad & \text{otherwise.}
  \end{cases}
\]
For example,
\[ \partial_{\omega_4}(\omega_2)=0; \qquad
   \partial_{\omega_2}(\omega_4)=
     8\mathcenter{%
\input{draws/SideGlu.tex}}
    +4\mathcenter{%
\input{draws/DiagGlu.tex}}.
\]
One might think of $D$ as a monomial of degree equal to the number
of legs. If $C$ has $k$ legs and degree $m$, then $\partial_{C}$
is an operator of degree $m-k$. By linear extension, we find that
every strutless  $C\in\cB$ defines an operator
$\partial_{C}\:\cB\to\cB$.  (We restrict to diagrams without
struts to avoid circles arising from the pairing of two struts
and to guarantee convergence: gluing with a strut lowers the
degree of a diagram, and so the pairing would not extend from
$\cB^f$ to $\cB$.)

In some sense, $\partial_ C$ is a diagrammatic analogue of
 a  constant coefficient differential operator. For
instance, one has:
\begin{itemize}
\item A diagram $C$ with $k$ legs reduces the number of legs by $k$,
      corresponding to a differential operator of order $k$.
\item If $k = 1$ ($C$ has only one leg), we have a Leibniz rule like
      that for linear differential operators:
\[
\partial_C(D_1 \sqcup D_2) = \partial_C(D_1) \sqcup D_2 + D_1 \sqcup\partial_C(D_2).
\]
      (Actually, all diagrams with only one leg are 0 in $\cB$, so
      we have to extend our space of diagrams slightly for this
      equation to be non-empty.  Adding some extra vertices of
      valence 1 satisfying no relations is sufficient.)
\item Multiplication on the differential operator side is the same
      thing as composition:
\begin{equation}\label{product}
\partial_{C_1 \sqcup C_2} = \partial_{C_1} \circ \partial_{C_2}.
\end{equation}
\end{itemize}

\section{Cabling}

\label{cabling}

The behaviour of cabling will be crucial to the proofs of all of
the Theorems of this paper.  In this section, we will review some
results of~\cite{LeMurakami:Parallel} on disconnected cabling and
prove a new result on connected cabling.

\subsection{Tangles, framed tangles, and the Kontsevich integral}

Suppose $X$ is a compact oriented 1-manifold. A {\em
tangle} with skeleton $X$ is a smooth proper embedding of $X$ into
$\RR\times \RR\times [0,1]\subset \RR^3$, considered up to isotopy
relative to the boundary. The Kontsevich of such a tangle takes
value in the space $\cA'(X)$, obtained from $\cA(X)$ by dividing
by the {\em framing independence} relation which says that a
diagram containing an isolated chord is equal to 0 (see
\cite{BarNatan:OnVassiliev}, we will not need $\cA'(X)$ in the
future). When $X$ does not have any circle component, there is a
canonical embedding from $\cA'(X)$ into $\cA(X)$, and the
Kontsevich integral can be considered valued in $\cA(X)$.

The {\em framed} Kontsevich integral of a {\em framed} tangle
with skeleton $X$ takes value in $\cA(X)$ (no framing
independence relation here). For technical reasons we will define a
{\em framed tangle} as a tangle: (a) with boundary lying on two
lines, the upper one $\RR\times\{0\}\times\{1\}$ and the lower one
$\RR\times\{0\}\times\{0\}$, and (b) equipped with a non-zero
normal vector field which is standard $(0,1,0)$ at every boundary
point. Framed tangles are considered up to isotopy as usual. In
$\RR^3$ the set of framing of each component can be canonically
identified with $\ZZ$. The framed Kontsevich integral of a framed
tangle $L$ is denoted by $Z(L)$. (For details, see
\cite{BarNatan:OnVassiliev, BarNatan2, LeMurakami:Universal,
LeMurakami:Parallel}. In \cite{LeMurakami:Universal,
LeMurakami:Parallel}, $Z(L)$ is denoted by $\hat Z_f(L)$.)

If a framed tangle $L'$ is obtained from another $L$ by
increasing the framing of a component labeled $x$ by 1, then we
have the following {\em framing formula}:
\begin{equation}\label{framing}
Z(L') = Z(L) \connect
\exp\left(\frac{1}{2}\mlarge\isolatedchord\right).
\end{equation}
where the connected sum is done on the component labeled $x$ and
$\isolatedchord\in \cA(\uparrow)\equiv \cA(\directedcircle)$ is
the Jacobi diagram based on $\directedcircle$ with one strut.

The framed Kontsevich integral depends on the positions of the
boundary points. To get rid of this dependence one has to choose
standard positions for the boundary points. It turns out that the
best ``positions" are in a limit, when all the boundary points go
to one fixed point. (One has to regularize the Kontsevich
integral in the limit.) In the limit one has to keep track of the
order in which the boundary points go to the fixed point. This
leads to the notion of {\em parenthesized framed tangle}, or {\em
q-tangle} in \cite{LeMurakami:Universal},  -- a framed tangle with
a non-associative structure on each of the two sequences of
boundary points on the upper and lower lines. For details, see
\cite{BarNatan2, LeMurakami:Universal}.

In all framed tangles in this paper, we assume that a
non-associative structure is fixed. In many cases, there is only
one non-associative structure, or the non-associative structure is
clear from the context.

\subsection{Coproduct and Sliding property}

 Let\qquad\qquad
$
\Delta^x_{x_1\dots x_n}\: \cA(\uparrow_x\cup X)\to
    \cA(\uparrow_{x_1}\cup\cdots\cup\uparrow_{x_n}\cup X)
$

or\qquad\qquad
$
\Delta^x_{x_1\dots x_n}\: \cA(\directedcircle_x\cup X)\to
    \cA(\directedcircle_{x_1}\cup\cdots\cup\directedcircle_{x_n}\cup\,X)
$

be the pull-back of the $n$-fold {\em disconnected} cover of the
component labeled $x$.  When we do not care about the labels on
the result, an alternate notation is $\Delta^{(n)}_x$.

Suppose $D\in
    \cA(\uparrow_{x_1}\cup\cdots\cup\uparrow_{x_n})$ and $D'\in
    \cA(\uparrow_{x_1}\cup\cdots\cup\uparrow_{x_n}\cup X)$. We
    define $D\cdot D'\in
    \cA(\uparrow_{x_1}\cup\cdots\cup\uparrow_{x_n}\cup X)$ as the
    element obtained by placing $D$ on top of $D'$, ie,
    identifying the lower endpoint of $\uparrow_{x_i}$ in $D$ with
 the upper endpoint of $\uparrow_{x_i}$ in $D'$, for $i=1,\dots, n$.
 Similarly, $D'\cdot D\in
    \cA(\uparrow_{x_1}\cup\cdots\cup\uparrow_{x_n}\cup X)$ is
    obtained by placing $D'$ on top of $D$.

    In general $D\cdot D' \neq D'\cdot D$. The following is a
    special case when one has equality (see, for example,
    \cite[Lemma 8.1]{LeMurakami:Parallel}):

    \begin{lemma}[Sliding property] The image of $\Delta_x^{(n)}$
    commutes with $
    \cA(\uparrow_{x_1}\cup\cdots\cup\uparrow_{x_n})$, ie, for
    every $D\in
    \cA(\uparrow_{x_1}\cup\cdots\cup\uparrow_{x_n})$ and $D' \in \cA(\uparrow_x\cup
    X)$, we have that $D \cdot \Delta_x^{(n)}(D') = \Delta_x^{(n)}
    (D') \cdot D$.
    \label{slid}
    \end{lemma}

With the above product,
$\cA(\uparrow_{x_1}\cup\cdots\cup\uparrow_{x_n})$ is an algebra.
There is also a co-product on
$\cA(\uparrow_{x_1}\cup\cdots\cup\uparrow_{x_n})$ which gives us
a structure of a Hopf algebra, and
$\cA(\uparrow_{x_1}\cup\cdots\cup\uparrow_{x_n})$ is a
(completed) polynomial algebra generated by primitive elements.
The isolated chord diagrams are among primitive elements. This is
the reason why there is a canonical algebra embedding from
$\cA'(\uparrow_{x_1}\cup\cdots\cup\uparrow_{x_n})$ into
$\cA(\uparrow_{x_1}\cup\cdots\cup\uparrow_{x_n})$.

\subsection{Disconnected cabling} Suppose $L$ is a framed tangle, with one of its components
labeled  $x$. The $n$-fold disconnected cabling of $L$ along $x$,
denoted by $\Delta_x^{(n)}(L)$, is  the tangle obtained from $L$
by replacing the component labeled $x$ with $n$ of its parallels.
Here the parallels are determined by the framing, and each
inherits a natural framing from that of component $x$.

The following proposition, proved in \cite{LeMurakami:Parallel},
describes the behaviour of the Kontsevich integral under
disconnected cabling.

\begin{proposition}
Suppose  that  a component labeled $x$ in a framed tangle $L$ is
either closed or has one upper and one lower boundary points.
Then
\begin{equation}\label{disconn-cable}
 Z(\Delta_x^{(n)} L) = \Delta_x^{(n)} (Z(L)).
\end{equation}
\end{proposition}

Since $Z(\Delta_x^{(n)} L)$ depends on the positions of the
boundary points, one has to be careful about the boundary points
of the new components (ie\ parallels) in $\Delta_x^{(n)} L$ when
the components label $x$ is not closed. The correct choice is the
one in which the distances between the boundary points of the
parallels are infinitesimally small compared to the distance
between any of these points and any other boundary point. In the
language of parenthesized framed tangles (or $q$-tangles), this
means the boundary points of the parallels must form an innermost
structure in the overall non-associative structure of the tangle
$L'$, and the non-associative structure among the boundary points
of the parallels on the upper line must be the same as that among
the boundary points of the parallels on the lower line.

\begin{remark} In general, the disconnected cabling formula
(\ref{disconn-cable}) does not hold true if the $x$ component has
both boundary points on the same upper or lower line. However, it
would hold true if the framed Kontsevich integral is modified by
using a good enough associator \cite{LeMurakami:Parallel}.
\end{remark}

\subsection{Connected cabling}

Let us define
\[
\psi^{(n)}_x\: \cA(\directedcircle_x\cup
X)\to\cA(\directedcircle_x\cup X)
\]
as the pull-back of the $n$-fold connected cover of the circle
labeled $x$.

Suppose $L$ is a framed tangle, with one of its {\em closed}
components labeled by $x$. The $n$-fold connected cabling of $L$
along $x$, denoted by $\CCable_x^{(n)}(L)$, is defined as follows.
On the torus boundary of a small tubular neighborhood of component
$x$ there are the preferred longitude and meridian. Replace the
component $x$ with a closed curve on the torus boundary whose
homology class is equal to that of the meridian plus $n$ times
the longitude. The result is $\CCable_x^{(n)}(L)$. The new
component inherits the orientation and framing from the old one.

The following theorem  describes the behaviour of the Kontsevich
integral under connected cabling.

\begin{theorem}
\label{conn-cable} Suppose a component labeled $x$ in a framed
tangle $L$ is closed (ie\ a knot).  Then
\[
Z(\CCable_x^{(n)}(L)) =
    \left[\psi^{(n)}_x(Z(L) \connect_x
    \exp(\frac{1}{2n}\isolatedchord))\right]\,
    \connect \exp(-\frac{1}{2}\isolatedchord).
\]
\end{theorem}

\proof
We will prove the theorem in the case when $L$ is a knot. The
case of an arbitrary tangle is quite similar.

 The difference
between the connected cabling and the disconnected cabling is the
extra $1/n$ twist $T_n$ inserted at one point:
\[
T_n = \figcent{1n-twist}{0.25}.
\]
By isotopy we can assume that this twist occurs in a horizontal
slice where all the other strands are vertical.  We can apply
(\ref{disconn-cable}) on the $(n,n)$ ``tangle'' obtained by
excising $T_n$.  (This object is not properly a tangle, since
there is a little piece cut out of it.  But we can still compute
its Kontsevich integral.) To complete the computation, we need to
compute $a := Z(T_n)$.

Repeating $T_n$ $n$ times, we get a full twist which we can
compute using the framing  and the disconnected cabling formulas
(\ref{framing}), (\ref{disconn-cable}):
\[
Z\left(\figcent{n-twists}{0.25}\right)
    = Z\left(\figcent{fulltw}{0.25}\right)
    = \Delta^{(n)}\left(\exp(\frac{1}{2}\isolatedchord)\right)
       \cdot \exp(-\frac{1}{2}\isolatedchord)^{\otimes n} =: b.
\]
The notation $\exp(-\frac{1}{2}\isolatedchord)^{\otimes n}$ means
$n$ copies of the framing change element
$\exp(-\frac{1}{2}\isolatedchord)$, one on on each of the $n$
strands, and the product $\cdot$ is the product in
$\cA(\uparrow\cup \dots \cup \uparrow)$.

The $n$ copies of $T_n$ that appear are not quite the same: they
differ by cyclic permutations of the strands. If we could arrange
the $n$ strands at the top and bottom of $T_n$ to be at the
vertices of a regular $n$-gon, the strands would be symmetric and
$a^n = b$ or
\[
a = b^{\frac1n} =
\Delta^{(n)}\left(\exp(\frac{1}{2n}\isolatedchord)\right)
        \exp(-\frac{1}{2n}\isolatedchord)^{\otimes n}.
\]
In reality, $a$ is not symmetric, ie, $\sigma(a) \neq a$, where
$\sigma$ is the automorphism of $\cA(\uparrow_{x_1}
\dots\uparrow_{x_n})$ which rotates the strands by $x_i \mapsto
x_{i-1}$. We have
$$Z\left(\figcent{n-twists}{0.25}\right) = a \cdot \sigma(a) \cdot
\sigma^2(a) \dots \sigma^{n-1}(a).$$
We can conjugate $T_n$ by some tangle $C$ to get the strands
symmetric: $T_n = C T'_n C^{-1}$, with $T_n'$ symmetric. From the
definition of the framed Kontsevich integral
\cite{LeMurakami:Universal}, it follows that $a = c \cdot a'
\cdot \sigma (c^{-1})$, where $a'= Z'(T'_n)$ is the usual
Kontsevich integral of $T'_n$, and $c= Z(C)\in \cA(\uparrow_{x_1}
\cup \dots \cup \uparrow_{x_n})$. Thus
$$ a \cdot \sigma(a) \cdot
\sigma^2(a) \dots \sigma^{n-1}(a) = c (a')^n c^{-1}.$$
And hence
\begin{align*}
a &= c \cdot b^{\frac{1}{n}} \cdot c^{-1}\\
& = c \cdot
\Delta^{(n)}\left(\exp(\frac{1}{2n}\isolatedchord)\right) \cdot
        \exp(-\frac{1}{2n}\isolatedchord)^{\otimes n} \cdot
c^{-1}.
\end{align*}
By the above
computations, the invariant of the connected cable of a knot $L$
is $\Delta^{(n)}(Z(L))$, multiplied by $a=Z(T_n)$, and close up
with a twist.  The conjugating elements $c$ and $c^{-1}$ can be
swept through the knot, using the sliding property of Lemma
\ref{slid}, and cancel each other. The factor
$\Delta^{(n)}(\exp(\isolatedchord/2n))$ in $a$ can be combined
with $Z(L)$ so that we apply $\Delta^{(n)}$ to $Z(L) \connect
\exp(\isolatedchord/2n)$.  The twisted closure turns
$\Delta^{(n)}$ into $\psi^{(n)}$.  The remaining $n$ factors of
$\exp(-\isolatedchord/2n)$ in $a$ can be slid around the knot and
combined to give
$$
Z(\CCable^n(L)) =
\left[\psi^{(n)}\left(Z(L)\connect\exp(\frac{1}{2n}\isolatedchord)\right)\right]
            \connect \exp(-\frac{1}{2}\isolatedchord).
\eqno{\qed}$$

\begin{remark} Suppose $\CCable^{(n|m)}_x(L)$ is the connected $(n,m)$-cabling of a framed tangle $L$
along a closed component labeled $x$, where $n$ and $m$ are
co-prime integer with $n>0$, ie, $\CCable^{(n|m)}_x(L)$ is
obtained by replacing the $x$ component with a closed curve on
the torus boundary of the regular neighborhood which represents
the homology class of $m$ times the meridian plus $n$ times the
longitude. Let $\psi^{(n|m)}$ denotes the corresponding pull-back
of Jacobi diagrams. Then the proof of Theorem \ref{conn-cable}
also gives:
\[
Z(\CCable^{(n|m)}(L)) = \left[
\psi^{(n|m)}\left(Z(L)\connect\exp(\frac{m}{2n}\isolatedchord)\right)\right]
            \connect \exp(-\frac{m}{2}\isolatedchord).
\]
\end{remark}

\subsection{Operators $\Delta,\psi$ and symmetrized diagrams}

One reason to introduce symmetrized diagrams is that the
operations $\Delta$ and $\psi$ above become very simple in $\cB$.
Using the symmetrization map one can trade an interval in the
skeleton with an asterisk, and a circle with a circled asterisk.
The map $\Delta_x$ and $\psi_x$ can be carried over to the new
spaces. The following lemmas are well-known  (and easy to check).

\begin{lemma}\label{lem:delta} The map
\[
\Delta^x_{x_1\dots x_n}\:\cA(\ast_x \cup X) \to
    \cA(\ast_{x_1} \cup\dots\cup\ast_{x_n} \cup X)
\]
is the sum over all ways of replacing each $x$ leg by one of the
$x_i$.\qed
\end{lemma}

\begin{remark}
$\Delta$ is similar to a coassociative, cocommutative coproduct
in a coalgebra, except that it does not take values in $\cA
\otimes \cA$.
\end{remark}

The operation $\Delta$ in Lemma~\ref{lem:delta} is analogous to a
change of variables $x \mapsto x_1 + \dots + x_n$ for ordinary
functions $f(x)$.  We will use a suggestive notation: a leg
labeled by a linear combination of variables means the sum over
all ways of picking a variable from the linear combination.  If
$D(x)$ is a diagram with some legs labeled $x$,
$\Delta^{(n)}(D(x)) = D(x_1 + \dots + x_n)$ is the diagram with
the same legs labeled $x_1 + \dots + x_n$.

\begin{lemma}[See \cite{Kricker}] \label{adams} The map
\[
\psi^{(n)}_x\: \cA(\circledast_x \cup X) \to
    \cA(\circledast_x \cup X)
\]
is multiplication by $n^k$ on diagrams with $k$ legs labeled $x$.
\qed
\end{lemma}

This operation is related to the change of variables $x \mapsto
nx$.

\section{The Wheeling Theorem}

\label{wheeling}

The operator  $\partial_\Omega\:\cB \to \cB$, where $\Omega$ is
the wheels element of the Introduction,  is called the
``wheeling'' map. The proof of the following theorem will occupy
the rest of this section.

    \begin{theorem}[Wheeling]\label{thm:wheeling}
The map $\Upsilon = \chi \circ \partial_\Omega\: \cB \to
\cA$ is an algebra isomorphism.  \end{theorem}

The map $\Upsilon$ is the diagrammatic analogue of the
Duflo--Kirillov map. Note that by (\ref{product}),
$\partial_\Omega \partial_{\Omega^{-1}}=\text{id}$, hence
$\partial_\Omega$ is a vector space isomorphism.  Since $\chi$ is
also a vector space isomorphism, $\Upsilon$ is automatically
bijective.

\subsection{An inner product}  Suppose $C, C' \in \cB$ are diagrams
such that $C$ has no struts. If  $C$ and $C'$ have the same number of
legs, then the {\em inner product} $\langle C, C' \rangle$ is the sum
of all ways of gluing all the legs of $C$ to all legs of $C'$. If $C$
and $C'$ do not have the same number of legs, then define $\langle C,
C' \rangle=0$. The restriction that $C$ not have struts is to guarantee
convergence and avoid closed circles.

       We will sometimes want to fix $C$ and consider $\langle
C\comma\cdot\rangle$ as a map from $\cB$ to $\cA(\emptyset)$; we
will denote this map $\iota(C)$.  This definition works equally
well in the presence of other skeleton components or to glue
several components.  We will use subscripts to indicate which
ends are glued.

There are two dualities relating $\langle\cdot\comma\cdot\rangle$
with other operations we have defined.

\begin{lemma}\label{lem:mult-comult-dual}
Multiplication and comultiplication in $\cB$ are dual in the
sense that
\[
\langle C, D_1 \sqcup D_2 \rangle =
    \langle \Delta_{xy} C, (D_1)_x \otimes (D_2)_y\rangle_{xy}.
\]
Similar statements hold in the presence of other ends.
\end{lemma}
\begin{proof} The glued diagrams are the same on the two sides; we
either combine the legs of $D_1$ and $D_2$ into one set and then
glue with $C$, or we split the legs of $C$ into two pieces which
are then glued with $D_1$ and $D_2$. (Note that there are no
combinatorial factors to worry about: in both cases, we take the
sum over all
possibilities.)  
\end{proof}

\begin{lemma}\label{lem:mult-diffop-dual} Multiplication by a diagram
$B \in\cB$ and applying $B$ as a diagrammatic differential
operator are adjoint in the sense that
\[
\langle A \sqcup B, C\rangle = \langle A, \partial_B(C) \rangle.
\]
\end{lemma}
\begin{proof}
As before, the diagrams are the same on both sides.
\end{proof}

\subsection{The map $\Phi$} \label{sec:phi} Let $\openhopf{x}{z}$
be the tangle in Figure~\ref{fig:One}, which is a bead (labeled
$x$ here) on a wire (labeled $z$ here). Its Kontsevich integral
$Z(\openhopf{x}{z})$ takes values in $\cA(\uparrow_z,
\directedcircle_x)$.
 Symmetrizing the legs attached to the
bead $x$ as explained in \ref{sec:diagrams2}, we get
\[
\chi^{-1}_x Z(\openhopf{x}{z}) \in \cA(\uparrow_z, \circledast_x).
\]
Finally, we use the inner product operation along the legs $x$ to
get a map from $\cB$ to $\cA$:
\[
\Phi = \iota_x \chi^{-1}_x Z(\openhopf{x}{z})\: \cB \to \cA
\]
In this last step, there are two things we have to check.  First,
we must see that $\chi^{-1}Z(\OpenHopf)$ has no struts. This
follows from the fact that we took the bead with the zero
framing.  Second, we need to check that the inner product
descends modulo the link relations on $x$ in
$\cA(\uparrow_z,\circledast_x)$.

\begin{lemma}\label{lem:glue-welldef}
The inner product
$\langle\cdot\comma\cdot\rangle_x\:\cA(\ast_x\cup X) \otimes
\cA(\ast_x)\to\cA(X)$ descends to a map
$\langle\cdot\comma\cdot\rangle_x\:\cA(\circledast_x\cup X)
\otimes \cA(\ast_x) \to \cA(X)$.
\end{lemma}

\begin{figure}[ht!]
\[
\mathcenter{\fig{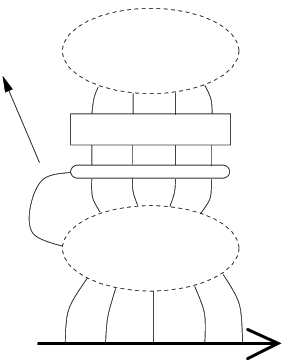}{0.28}}=
\mathcenter{\fig{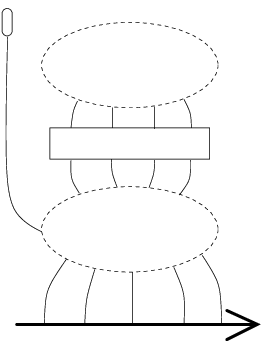}{0.28}}= 0.
\]
\caption{The proof that $\Phi(D)$ is well-defined modulo link
relations on $Z(\OpenHopf)$: link relations in $Z(\OpenHopf)$ can
be
slid over $D$.}%
\label{fig:phi-welldef}
\end{figure}

\begin{proof}  Link relations in
$\cA(\circledast_x\sqcup X)$ can be slid over diagrams in
$\cA(\ast_x)$, as shown in Figure~\ref{fig:phi-welldef}. (See
similar arguments in \cite{BarNatan:OnVassiliev}).
\end{proof}

\subsection{Multiplicativity of $\Phi$}\label{sec:multiplicativity}
We now come to the key lemma in the proof of the wheeling theorem.

\begin{lemma}\label{lem:phi-multiplicative} The map $\Phi\: \cB
\to \cA$ is an algebra map.
\end{lemma}

\begin{figure}[ht!]
\begin{center}
\fig{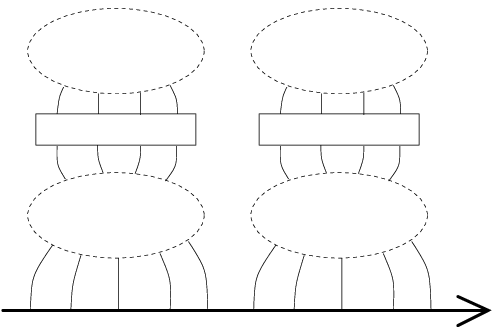}{0.28}
\end{center}
\caption{Gluing $Z(\OpenHopf\connect\OpenHopf)$ to $D_1 \otimes
           D_2$}\label{fig:glue-1+1}%
\end{figure}

\proof
As advertised, we use the equality of links ``$1+1=2$''.  Let us
see what this equality of links says about the Kontsevich integral
of the Hopf link.  On the ``$1+1$'' side, we see the connected
sum of two open Hopf links. It is known that the invariant of the
connected sum is the connected sum of the invariants.  To write
this conveniently, let $H(z;x)$ be $Z(\OpenHopf)
\in\cA(\uparrow_z, \circledast_x)$, with the wire labeled by $z$
and the bead labeled by $x$.  Then
\[
Z(\OpenHopf \connect \OpenHopf) = H(z; x_1) \connect_z H(z; x_2)
  \in \cA(\uparrow_z, \circledast_{x_1}, \circledast_{x_2}).
\]
On the ``$2$'' side, we see the disconnected cable of a Hopf
link.  By the disconnected cabling formula (\ref{disconn-cable}),
this becomes the coproduct $\Delta$:
\[
Z(\Delta_x^{(2)}(\OpenHopf)) = \Delta^x_{x_1x_2}H(z;x)
  \in \cA(\uparrow_z, \circledast_{x_1}, \circledast_{x_2}).
\]
Since the two tangles are isotopic, we have
\begin{equation} \label{eq:HEquation}
H(z; x_1, x_2) \overset{\text{def}}{=}
    H(z;x_1) \connect_z H(z;x_2) = \Delta^x_{x_1x_2}H(z;x)
  \in \cA(\uparrow_z, \circledast_{x_1}, \circledast_{x_2}).
\end{equation}
Now consider the map
\[
\Xi = \iota_{x_1}\iota_{x_2}H(z; x_1, x_2)\: \cB \otimes \cB
\to \cA;
\]
in other words, in $\Xi(D_1 \otimes D_2)$ glue the $x_1$ and $x_2$
legs of $H(z;x_1,x_2)$ to $D_1$ and $D_2$ respectively.  This
descends modulo the two different link relations in
$\cA(\uparrow, \circledast, \circledast)$ by the argument of
Figure~\ref{fig:phi-welldef}, applied to $D_1$ and $D_2$
separately.  We have two different expressions for this map from
the two different expressions for $H(z; x_1, x_2)$.  On the
``$1+1$'' side, the gluing does not interact with the connected
sum and we have
\[
\Xi(D_1, D_2) = \Phi(D_1) \connect \Phi(D_2),
\]
see Figure~\ref{fig:glue-1+1}. 

\begin{figure}[ht!]
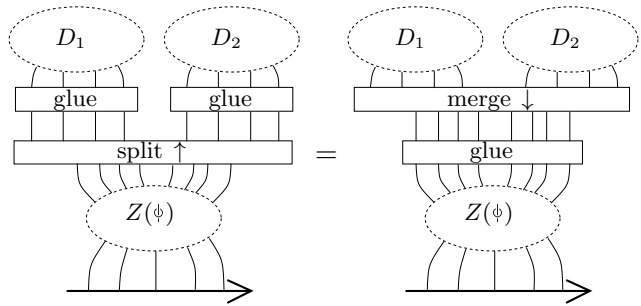

\begin{center}
\[
\figcent{glue-2a}{0.28} = \figcent{glue-2b}{0.28}
\]
\end{center}
\caption{Gluing $Z(\Delta_x^{(2)}(\OpenHopf))$ to $D_1 \otimes D_2$ in two
equivalent ways}%
    \label{fig:glue-2}%
\end{figure}

For the ``$2$'' side, we use
Lemma~\ref{lem:mult-comult-dual} to see that
\begin{align}
\Xi(D_1, D_2) &= \langle \Delta^x H(z; x), D_1 \otimes D_2 \rangle \\
          &= \langle H(z; x), D_1 \sqcup D_2 \rangle \\
          &= \Phi(D_1 \sqcup D_2),
\end{align}
see Figure~\ref{fig:glue-2}.

Combining the two, we find
$$
\Phi(D_1) \connect \Phi(D_2) = \Phi(D_1 \sqcup D_2).
\eqno{\qed}$$

\subsection{Mapping degrees and the Duflo--Kirillov
isomorphism}\label{sec:phi0} We have successfully constructed a
multiplicative map from $\cB$ to $\cA$.  We will see later  that
this map $\Phi$ is the same as $\Upsilon$, but we cannot yet see
this. Instead we will consider the lowest degree term $\Phi_0$ of
$\Phi$.

The {\em mapping degree} of a diagram $D \in \cA(\uparrow_z,
\circledast_x)$ with respect to $x$ is the amount $\iota_x D\: \cB
\to \cA$ shifts the degree.  Explicitly, it is the degree
of~$D$ minus the number of $x$ legs of~$D$.

Since there are no $x$-$x$ struts in $H(z;x)$, every $x$ leg of
$H$ must be attached to another vertex (either internal or on the
interval $z$).  Furthermore, if two $x$ legs are attached to the
same internal vertex, the diagram vanishes by anti-symmetry.
Therefore there are at least as many other vertices as $x$ legs
in $H$ and the mapping degree is $\geq 0$.

Let $H_0(z;x)$ be the part of $H(z;x)$ of mapping degree 0 with
respect to $x$, and   $\Phi_0\: \cB \to \cA$ be $\iota_x H_0(z;x)$.
The map $\Phi_0$ is still multiplicative, since the
multiplications in $\cA$ and $\cB$ both preserve degrees.  (For
homogeneous diagrams $D_1$ and $D_2$ of degrees $n_1$ and $n_2$,
$\Phi_0(D_1 \sqcup D_2)$ is the piece of $\Phi(D_1 \sqcup D_2)$
of degree $n_1+n_2$ and likewise for $\Phi_0(D_1) \connect
\Phi_0(D_2)$.)

\begin{figure}[ht!]
\begin{center}
\fig{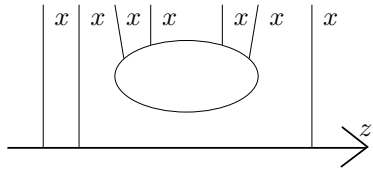}{0.75}
\end{center}\vspace{-20pt}
\caption{The only diagrams in $\cA(\uparrow_z, \circledast_x)$ of
mapping degree 0 with respect to $x$ are wheels and struts.}
\label{fig:wheels-struts}
\end{figure}

The diagrams that appear in $H_0$ are very restricted, since every
vertex that is not an $x$ leg must connect to an $x$ leg.  The
possible diagrams are $x$ wheels and $x-z$ struts, as shown in
Figure~\ref{fig:wheels-struts}. The linking number between the
bead and the wire in the link $\OpenHopf$ is 1, so the
coefficient of the $x$-$z$ strut is 1. Combined with the fact
that the Kontsevich integral is grouplike
\cite{LeMurakami:Parallel}, we find that
\begin{align*}
H_0(z;x) &= \exp(\strutv{x}{z}) \sqcup \Omega', \quad \text{where} \\
\Omega' &= \exp_\sqcup (\sum_n a_{2n}\omega_{2n})
\end{align*}
for some coefficients $a_{2n}$.  Note that the right hand side is
written in $\cA(\uparrow, \ast)$ (with a strange mixed product),
since there is no algebra structure on $\cA(\uparrow,
\circledast)$.

By the following lemma, we now have a multiplicative map very
similar to our desired map~$\Upsilon$.

\begin{lemma} One has that
$\Phi_0 = \chi \circ \partial_{\Omega'}$.
\end{lemma}

\proof
Using Lemma~\ref{lem:mult-diffop-dual} and noting that gluing
with $\exp(\strutv{x}{z})$ takes the legs of a diagram in $\cB$
and averages over all ways of ordering them, as in the definition
of $\chi$, we see that
$$
\Phi_0(D) = \langle \exp(\strutv{x}{z} \sqcup \Omega', D \rangle
      = \langle \exp(\strutv{x}{z}, \partial_\Omega'(D) \rangle
      = \chi(\partial_\Omega'(D)).\eqno{\qed}$$

\subsection{Identifying $\Phi_0$ with $\Upsilon$} To complete the
proof of the wheeling theorem, one needs only to show that
$\Omega=\Omega'$, or $a_n=b_n$ for $n=2,4,6,\dots$. This can be
proved as follows.

First of all, a calculation of the degree 2 part of the Kontsevich
integral of the Hopf link will show that $a_2=b_2$. Thus if
$\Omega \neq \Omega'$, then for some $n>1$,
$$\Omega^{-1} \Omega' = 1 + (a_{2n}-b_{2n})\omega_{2n}+ \text{higher order terms}.$$
Second, the map $\Upsilon = \chi \circ \partial_\Omega$ is known
to be an algebra isomorphism on the level of simple Lie algebras
\cite{BGRT:WheelsWheeling}. Thus for a simple Lie algebra
$\frakg$, the map $\partial_{\Omega^{-1}\Omega'}$ is an algebra
automorphism of $S(\frakg)^{\frakg}$. When $\frakg=sl_2$, the
algebra $S(\frakg)^{\frakg}$ is a polynomial algebra on one
generator, which is the image of the strut $\strutn{}{}$. On the
strut $\partial_{\Omega^{-1}\Omega'}$ acts as the identity (since
there is no non-trivial diagram with less than 3 legs), hence
$\partial_{\Omega^{-1}\Omega'}$ acts as the identity on the whole
algebra $S(\frakg)^{\frakg}$.

Third, the action of $\omega_{2n}$ on $S(\frakg)^{\frakg}$ is
non-trivial. Explicitly, $\partial_{\omega_{2n}}[(\strutn{}{})^n ]
= 2 (2n+1)!$ in $sl_2$, which can be proved easily by induction.
Thus, if $a_{2n} \neq b_{2n}$, then
$\partial_{\Omega^{-1}\Omega'}$ cannot act as identity on
$S(\frakg)^{\frakg}$.

We conclude that $\Omega= \Omega'$, and this completes the proof
of the wheeling theorem. For another proof of $\Omega= \Omega'$,
more detailed and without using the result of
\cite{BGRT:WheelsWheeling}, see the Appendix.

\subsection{Back to the Duflo--Kirillov isomorphism} We note that the
wheeling theorem implies the multiplicative property of the
Duflo--Kirillov isomorphism for a metrized Lie (super-) algebra
$\frakg$.  Indeed, using the standard maps $W_\frakg$ from spaces of
diagrams into spaces of tensors, we set $J=W_\frakg(H(z; x))\in
U(\frakg)\otimes S(\frakg)_\frakg$. Here $S(\frakg)_\frakg$ denotes the
space of coinvariants of the $\frakg$ action on $S(\frakg)$ --- the
link relation dictates the descent to this quotient of $S(\frakg)$.
Also, strictly speaking $J$ lives in the completion of
$U(\frakg)\otimes S(\frakg)_\frakg$ induced by the grading on
$S(\frakg)_\frakg$.  Equation~(\ref{eq:HEquation}) and the
compatibility between $W_\frakg$ and multiplication and
comultiplication imply now that $J$ satisfies
\begin{equation} \label{eq:JEquation}
  J\connect J = (1\otimes\Delta)J
  \qquad\text{in}\qquad
  U(\frakg)\otimes S(\frakg)_\frakg\otimes S(\frakg)_\frakg,
\end{equation}
where $J\connect J$ denotes the result of multiplying two copies of $J$
using the product of $U(\frakg)$, so that the result is in (the
appropriate completion of) $U(\frakg)\otimes S(\frakg)_\frakg\otimes
S(\frakg)_\frakg$. Now use the metric of $\frakg$ to identify the space
of coinvariants in $S(\frakg)$ as the dual of the space
$S(\frakg)^\frakg$ of invariants and hence to re-interpret $J$ as an
element of $U(\frakg)\otimes(S(\frakg)^\frakg)^\star$ and hence as a
map $W_J\:S(\frakg)^\frakg\to U(\frakg)$. One easily verifies that
equation~(\ref{eq:JEquation}) implies that $W_J$ is multiplicative. It
remains to see that $W_J$ is equal to the Duflo--Kirillov map $\chi
\circ \partial_{j^{\frac12}}$. This follows from the computation of
$H(z; x)$ in terms of the diagrammatic analogue $\Omega$ of
$j^{\frac12}$ in Section~\ref{Hopf}.

\section{The Wheels Theorem. The Kontsevich integral of the unknot}\label{wheels}

This section is devoted to the proof of the Wheels theorem.
\begin{theorem}[Wheels]
The framed Kontsevich integral of the unknot is the wheels
element:
\[
Z(\directedcircle) = \chi (\Omega).
\]
\end{theorem}

 We will denote $\nu =
Z(\directedcircle)\in \cA = \cA(\uparrow)\equiv
\cA(\directedcircle)$.

\subsection{Useful facts}

We will first derive some nice properties the wheels element
$\Omega$. Set $H_0(z;x) = \Omega_x \exp(\strutn{x}{z})$ and start from
the basic equality proved in the Wheeling theorem,
\[
\Delta^x_{x_1x_2}H_0(z;x) = H_0(z;x_1) \connect_z H_0(z;x_2)
  \in \cA(\uparrow_z \circledast_{x_1} \circledast_{x_2}).
\]
Now consider dropping the strand $z$, ie, mapping all diagrams
with a $z$ vertex to 0.  (Knot-theoretically, this corresponds to
dropping the central strand in the equation ``$1+1=2$''.)  We find
\begin{equation}\label{eq:delta-omega}
\Delta \Omega = \Omega \tensor \Omega \in \cA(\circledast
\circledast).
\end{equation}
Note that this equality is not true inside $\cA(\ast\ast)$.

\begin{lemma}[Pseudo-linearity of $\log\Omega$, see also \cite{Bar-NatanLawrence:RationalSurgery}]\label{lem:pseudo-linear}
For any $D \in \cB$,
\[
\partial_D(\Omega)
    = \langle D, \Omega \rangle \Omega.
\]
\end{lemma}
\begin{proof}
\[ \partial_D(\Omega)_x = \langle D_y, \Omega_{x+y}\rangle_y
                = \langle D_y, \Omega_x \Omega_y\rangle_y
                = \langle D_y, \Omega_y\rangle_y \Omega_x.
\]
In the second equality, we use Equation~(\ref{eq:delta-omega}).
This is allowed, since the contraction descends to
$\cA(\circledast \circledast) \isom \cA(\circledast \uparrow)$ by
the argument of Lemma~\ref{lem:glue-welldef}.
\end{proof}

\begin{remark} Compare this lemma with standard calculus: if $D$ is
any differential operator and $f$ is a linear function, then
$De^f=(Df)(0)e^f$.  The prefix ``pseudo'' is written above because
Lemma~\ref{lem:pseudo-linear} does not hold for every $D$, but
only for $x$-invariant $D$'s, ie, for $D$ with link relations on
$x$-legs.
\end{remark}

Although we are interested in knots and links in $S^3$, for which
the appropriate space of diagrams is the boundary connected part
$\cA^\bc$, vacuum diagrams (elements of $\cA(\emptyset)$) appear
at various points. Notably, the wheeling map $\Upsilon$ does not
preserve the subspace of boundary connected diagrams.
Although the resulting vacuum
components can be computed explicitly,%
\footnote{D.~Bar-Natan and
R.~Lawrence~\cite{Bar-NatanLawrence:RationalSurgery} have done
these computations} they are almost always irrelevant for us and
it would just complicate the formulas to keep track of them.  To
avoid this, we will introduce the {\em boundary-connected
projection} $\pi^\bc\: \cA \to \cA^\bc$ which maps any
diagram containing vacuum components to 0 and is otherwise the
identity.  Note that $\pi^\bc$ is multiplicative.  There are
similar projections, which will also be called $\pi^\bc$, for
other spaces $\cA(X)$.

If we compose Lemma~\ref{lem:pseudo-linear} with $\pi^\bc$, we
find, for a diagram $D\in \cB$,
\begin{equation}\label{eq:pseudo-linear-2}
\pi^\bc \partial_D \Omega = \begin{cases}
    \Omega  & \text{$D$ is the empty diagram} \\
    0   & \text{otherwise.}
  \end{cases}
\end{equation}

\subsection{A lemma on the bound of numbers of legs}
\begin{lemma}\label{lem:wheels-lemma} For any elements $x_1, \dots,
x_k \in \cA(\uparrow)$ with at least one leg on the interval
$\uparrow$, $\chi^{-1}(x_1 \connect \cdots \connect x_k) \in \cB$
has at least $k$ legs.
\end{lemma}

\begin{proof} First note that any vacuum diagrams that appear in the
$x_i$'s pass through unchanged to the result; let us assume that
there are none, so that we can use the vacuum projection
$\pi^\bc$ without changing the result.  By the wheeling theorem,
\[
\pi^\bc \chi^{-1}(x_1 \connect \dots \connect x_n)
    = \pi^\bc \partial_\Omega(\Upsilon^{-1}(x_1)\sqcup\dots\sqcup\Upsilon^{-1}(x_k)).
\]
Let $y_i = \pi^\bc \Upsilon^{-1}(x_i)$.  Each $y_i$ has at least
one leg, since if the $\partial_\Omega^{-1}$ of $\Upsilon^{-1}=
\partial_\Omega^{-1} \chi^{-1}$ eats all the legs of $\chi^{-1}x_i$, it
also creates a vacuum diagram which is killed by $\pi^\bc$.  Then
\[
\pi^\bc \partial_\Omega(y_1 \dots y_k) =
    \pi^\bc \langle \Omega_a, \Delta_{ab}(y_1 \dots y_k)\rangle_a.
\]
Let $\Delta_{ab} y_i = (y_i)_a + z_i$; diagrams in $z_1$ have at
least one $b$ leg.  We see that
\begin{align*}
\pi^\bc\langle\Omega_a, (y_1)_a \Delta_{ab}(y_2 \dots
y_n)\rangle_a
 &= \pi^\bc\langle(\partial_{y_1}\Omega)_a, \Delta_{ab}(y_2\dots y_k)\rangle_a
    &&\text{by Lemma~\ref{lem:mult-diffop-dual}} \\
 &= 0   &&\text{by Equation~\ref{eq:pseudo-linear-2}.}
\end{align*}
Therefore
\begin{align*}
\pi^\bc \partial_\Omega(y_1 \dots y_k)
   &= \pi^\bc(\langle \Omega_a, (y_1)_a \Delta_{ab}(y_2 \dots y_k)\rangle_a
    + \langle \Omega_a, z_1 \Delta_{ab}(y_2 \dots y_k)\rangle_a) \\
   &= \pi^\bc\langle \Omega_a, z_1 \Delta_{ab}(y_2 \dots y_k)\rangle_a \\
   & = \cdots \\
   & = \pi^\bc\langle \Omega_a, z_1 \dots z_k \rangle_a.
\end{align*}
Each $z_i$ has at least one leg labeled $b$, so the product has
at least $k$ legs labeled $b$ which are the legs in the result.
\end{proof}

\subsection{Coiling the unknot. Proof of the Wheels theorem}\label{sec:coiling}

The basic equation we will use to identify
$\nu=Z(\directedcircle)$ is ``$n\cdot 0 = 0$'' from the
introduction: the $n$-fold connected cable of the unknot is the
unknot with a new framing. The connected cabling formula of
Theorem \ref{conn-cable} implies that
\begin{equation}\label{eq:n00}
\psi^{(n)}(\nu \connect \exp_{\connect}(\frac{1}{2n}
\isolatedchord))
     = \nu \connect \exp_{\connect}(\frac{n}{2}\isolatedchord).
\end{equation}
This equation is true for all $n \in \ZZ$, $n > 0$.  In each
degree, each side is a Laurent polynomial in $n$ of bounded
degree; therefore, the two sides are equal as Laurent
polynomials.  The RHS is a polynomial in $n$, so both sides are
polynomials (ie, have no negative powers of $n$.)  Let us
evaluate both sides at $n=0$.  On the RHS, we get just $\nu$. For
the LHS, recall how $\psi^{(n)}$ acts in the space $\cB$: it
multiplies a diagram with $k$ legs by $n^k$ (see Lemma
\ref{adams}).

Consider expanding the exponential
$\exp_{\connect}(\isolatedchord/2n)$ in the LHS of
Equation~\ref{eq:n00}.  In the term with $(\isolatedchord)^k$,
there is a factor of $1/n^k$ from the coefficient $1/2n$.  On the
other hand, by Lemma~\ref{lem:wheels-lemma}, the product has at
least $k$ legs, or $k+1$ if there is a non-trivial contribution
from $\nu$. Since the overall power of $n$ is $n^{\text{\#
legs}-k}$, when we evaluate at $n=0$ the term $\nu$ does not
contribute at all. Hence
\[
\psi^{(n)}(\nu \connect
\exp_{\connect}(\frac{1}{2n}\isolatedchord))|_{n=0}
    = \psi^{(n)}(\exp_{\connect}(\frac{1}{2n}\isolatedchord))|_{n=0}.
\]
Now we want to pick out the term from $(\isolatedchord)^{\connect
k}$ with exactly $k$ legs.  We can do this computation explicitly
using the wheeling map $\Upsilon$. Alternatively, the result
 must be a diagram of degree $k$ and with $k$
legs, hence $\nu = \nu_0$, the part of mapping degree 0. It was
shown in Section~\ref{wheeling} that the  part of of mapping
degree 0 of  $Z(\OpenHopf)$ is $\Omega_x
\sqcup\exp_\sqcup(|_x^z)$. Dropping the central strand from
$\OpenHopf$ leaves an unknot, so $\Omega =\nu_0 = \nu$. This
completes the proof of the Wheels theorem.\qed

\begin{exercise} Do the computation suggested above.  Show that
\[
\chi^{-1}(\exp_{\connect}(\frac{1}{2}\isolatedchord)) =
    \Omega \sqcup \exp_{\sqcup} (\frac{1}{2}\strutn{}{}).
\]
\end{exercise}
\begin{hint} Use Lemma~\ref{lem:inspired-guess}.
\end{hint}

\section {From the unknot to the Hopf link}\label{Hopf}

By changing the framing on the unknot and cabling it, we can
construct a Hopf link.  Using the results of Section~\ref{cabling}
and the value of $Z(\bigcirc)$, we can compute the invariant of
the Hopf link from the invariant of the unknot. There are several
good formulas for the answer. An alternative exposition of the
results of this section can be found in
\cite{Bar-NatanLawrence:RationalSurgery}.

\begin{theorem}\label{thm:hopf-link}
The framed Kontsevich integral of the Hopf link can be expressed
in the following equivalent ways:
\begin{align*} Z(\hopf{x}{y}) &=
  \begin{cases}
    \Upsilon_x\circ\Upsilon_y(\exp(\strutn{y}{x})))\cdot \Vac \\
    \Upsilon_x(\exp_{\sqcup}(\strutn{y}{x})\Omega_x) \cdot \Vac
  \end{cases} \\
Z(\openhopf{x}{y}) &= \exp(\strutn{y}{x})\sqcup\Omega_y,
\end{align*} for some elements $\Vac \in \cA(\emptyset)$.
\end{theorem}

In the last expression, $\openhopf{x}{y}$ is the $(1,1)$ tangle
whose closure is the Hopf link, with the bead labeled by $y$ and
the wire labeled by $x$. From this last equality in
Theorem~\ref{thm:hopf-link}, we can see exactly the map $\Phi$
from Section~\ref{wheeling}.

\begin{corollary}
$\Phi = \Phi_0 = \chi \circ \partial_\Omega.$
\end{corollary}

\begin{proof} (of Theorem \ref{thm:hopf-link}) We start by computing the Kontsevich integral of the $+1$
framed unknot. In what follows we identify $\cB$ and $\cA$ using
$\chi$, and use $\sqcup$ and $\connect$ to denote the two
different products on $\cB$.
\begin{align*}
Z(\bigcircle^{+1}) &= \nu \connect \exp_{\connect}(\frac 12\isolatedchord)\\
    &= \partial_\Omega\left(\partial_\Omega^{-1}(\Omega) \sqcup
        \exp_{\sqcup}(\partial_\Omega^{-1}(\strutn{}{}))\right)
        && \text{by Theorem~\ref{thm:wheeling}} \\
    &= \pi^\bc\partial_\Omega\left(\Omega \sqcup \exp(\strutn{}{})\right).
        && \text{by Equation~\ref{eq:pseudo-linear-2}}
\end{align*}
To pass to the Hopf link, we double $Z(\bigcircle^{+1})$.  The
following lemma, which is obvious from the definition, tells us
how $\partial_\Omega$ interacts with doubling.  We use $\hat D$
as an alternate notation for $\partial_D$ so that we can use a
subscript to indicate which variable the differential operator
acts on.
\begin{lemma}\label{lem:double-diffop} For $C, D \in \cB$ with $C$ strutless,
\[
\Delta_{xy} \hat C(D) = \hat C_x(\Delta_{xy} D) = \hat
C_y(\Delta_{xy} D).
\]
\end{lemma}
If we want to apply $\partial_\Omega^{-1}$ to both components of
the Hopf link, we can compute
$\partial_\Omega^{-2}(Z(\bigcircle^{+1}))$.

\begin{lemma}\label{lem:inspired-guess}
 $\pi^\bc\partial_\Omega(\exp\frac{1}{2}\strutn{}{}) =
    \Omega\sqcup\exp(\frac{1}{2}\strutn{}{})$.
\end{lemma}
\proof
{\def\neg{\hspace{-8mm}}
\begin{align*}
\pi^\bc\partial_\Omega(\exp(\frac12\strutn{}{}))
   &= \pi^\bc \langle\Omega_y, \exp(\frac{1}{2} \strutn{x+y}{x+y})\rangle_y
\\
   &\neg= \pi^\bc\langle\Omega_y,
    \exp(\frac12\strutn xx)\exp(\strutn yx)\exp(\frac12\strutn yy)\rangle_y \\
   &\neg= \pi^\bc\langle \partial_{\exp(\frac12\strutn{}{})}(\Omega)_y,
    \exp(\strutn xy)\rangle_y\sqcup\exp(\frac12\strutn xx)
            && \text{by Lemma~\ref{lem:mult-diffop-dual}} \\
   &\neg= \pi^\bc\langle\Omega_y, \exp(\strutn yx)\rangle_y\sqcup\exp(\frac12\strutn xx)
            && \text{by Equation~\ref{eq:pseudo-linear-2}} \\
   &\neg= \Omega\sqcup\exp(\frac12\strutn{}{}).\tag*{\qed}
\end{align*}}

As a corollary, we see that
\begin{equation} \label{eq:inspired-guess}
\pi^\bc\partial_\Omega^{-2}(Z(\bigcircle^{+1})) =
\exp(\frac12\strutn{}{}).
\end{equation}
We now compute.
{\def\neg{\hspace{-0.2in}}
\begin{align*}
\pi^\bc \Delta_{xy}(\hat\Omega^{-2}Z(\bigcircle^{+2}))
   &= \pi^\bc \hat\Omega_x^{-1}\hat\Omega_y^{-1} Z(\,{}_x^{+1}\!\Hopf_y^{+1})
    &&\parbox{1.2in}{by Lemma~\ref{lem:double-diffop} and formula~(\ref{disconn-cable})} \\
   &\neg= \pi^\bc\Delta_{xy}(\exp(\frac12\strutn{}{}))
    &&\text{by Equation~\ref{eq:inspired-guess}} \\
   &\neg= \exp(\strutn{x}{y})\exp(\frac12\strutn xx)\exp(\frac12\strutn yy).
\end{align*}}
Apply $\Upsilon_x \circ \Upsilon_y$ to both sides.  We see that
\begin{align*}
Z(\,{}_x^{+1}\!\Hopf_y^{+1}) &= \pi^\bc Z(\,{}_x^{+1}\!\Hopf_y^{+1}) \\
    &= \pi^\bc\Upsilon_x \circ \Upsilon_y(\exp(\strutn{x}{y})\sqcup\exp(\frac12\strutn xx)
            \sqcup\exp(\frac12\strutn yy)) \\
    &= \pi^\bc \Upsilon_x \circ \Upsilon_y(\exp(\strutn{x}{y}))
        \connect \exp_{\connect}(\frac12\isolatedchord_x)
        \connect \exp_{\connect}(\frac12\isolatedchord_y) \\
\intertext{so} Z(\hopf{x}{y}) &= \pi^\bc \Upsilon_x \circ
\Upsilon_y(\exp(\strutn{x}{y})).
\end{align*}
This is the first equality of Theorem~\ref{thm:hopf-link}. For
the second equality,
\[
\hat\Omega_y(\exp(\strutn xy)) = \Omega_x \sqcup \exp(\strutn xy).
\]
so
\[
Z(\hopf{x}{y}) = \pi^\bc \Upsilon_x(\exp(\strutn{y}{x})\Omega_x).
\]
For the last equality of the theorem, multiplicativity of
$\Upsilon$ implies that
\begin{align*}
\pi^\bc\Upsilon_x(\exp(\strutn yx)\Omega_x)
  &= \pi^\bc(\Upsilon_x(\exp(\strutn yx)) \connect \Upsilon_x(\Omega_x)) \\
  &= \pi^\bc(\Upsilon_x(\exp(\strutn yx))) \connect \chi(\Omega_x) \\
  &= \chi(\exp(\strutn yx) \sqcup \Omega_y) \connect \chi(\Omega_x).
\end{align*}
Hence we have
\[
Z(\openhopf{x}{y}) = Z(\hopf{x}{y}) \connect \Omega_x^{-1}
    = \exp(\strutn yx) \sqcup \Omega_y.
\]
This completes the proof of Theorem~\ref{thm:hopf-link}.
\end{proof}

\section*{Appendix}

To show that $\Omega' = \Omega$, one can use the following
``Sawon's identity~\cite{HitchinSawon:CurvatureHyperKahler}'':
\begin{equation}
\label{eq:sawon} \langle \Omega', (\strutn{}{})^n \rangle =
(\frac{1}{24}\ThetaGraph)^n.
\end{equation}

\proof
Proceed by induction on $n$.  The result is trivial for $n=0$.
\begin{align*}
\langle \Omega', (\strutn{}{})^n\rangle
   &= \langle \Omega', \strutn{}{} \sqcup (\strutn{}{})^{n-1}\rangle \\
   &= \langle \partial_{\strutu{}{}} (\Omega'), (\strutn{}{})^{n-1}\rangle
    &&\text{by Lemma~\ref{lem:mult-diffop-dual}} \\
   &= \frac{1}{24}\ThetaGraph\,\langle\Omega', (\strutn{}{})^{n-1}\rangle
    &&\text{by Lemma~\ref{lem:pseudo-linear} and explicit computation} \\
   &= \left(\frac{1}{24}\ThetaGraph\right)^n
    &&\text{by induction}\tag*{\qed}
\end{align*}

The following is well-known, see eg \cite{Chmutov}.

\begin{lemma}\label{lem:sl2rel}
In the Lie algebra $\sltwo$, with the invariant inner product
$\langle x,y\rangle = -\tr(xy)$, where the trace is taken in the
adjoint representation, we have the following relations:
\begin{equation*}
\bigcirc \equiv 3 \qquad \IGraph \equiv \smoothing - \crossing
\end{equation*}
\end{lemma}

For example, apply the $\sltwo$ relations, we find that that
$\ThetaGraph \equiv 6$.

\begin{lemma}\label{lem:sl2-wheel}
Modulo the $\sltwo$ relations, $\omega_{2n} \equiv
2(\strutn{}{})^n$.
\end{lemma}
\proof
Proceed by induction.  This is a straightforward computation for
$n=1$.  For $n>1$, compute as follows:
$$\omega_{2n} = \figcent{sl2-w1}{0.12}
      = \figcent{sl2-w2}{0.12} - \figcent{sl2-w3}{0.12}
      = \strutn{}{} \sqcup \figcent{sl2-w4}{0.12}
      = \strutn{}{} \sqcup\omega_{2n-2}.\eqno{\qed}$$

\begin{lemma}\label{lem:sl2-gluing}
Modulo the $\sltwo$ relations, $\langle (\strutn{}{})^n,
(\strutn{}{})^n \rangle = (2n+1)!$.
\end{lemma}
\proof
Proceed by induction.  The statement is trivial for $n=0$.  For
$n>0$, the two ends of the first strut on the left hand side can
either connect to the two ends of a single right hand strut or
they can connect to two different struts.  These happen in $2n$
and $2n\cdot(2n-2)$ ways, respectively.  (Note that there are
$2n\cdot(2n-1)$ ways in all of gluing these two legs.)  We
therefore have
\[
\figcent{sl2-1}{0.15} = 2n\cdot\figcent{sl2-2}{0.15}
    + 2n \cdot (2n-2) \cdot \figcent{sl2-3}{0.15}
\]
and
\begin{align*}
\bigl\langle(\Strutn)^n, (\Strutn)^n\bigr\rangle
  &= \bigl(2n\bigcircle + 2n\cdot(2n-2)\bigr)
    \bigl\langle(\Strutn^{n-1},(\Strutn)^{n-1}\bigr\rangle \\
  &\equiv 2n\cdot(2n+1)\bigl\langle(\Strutn^{n-1},(\Strutn)^{n-1}\bigr\rangle \\
  &\equiv (2n+1)! &&\text{by induction.}\tag*{\qed}
\end{align*}

\begin{proposition}\label{prop:coefficients} One has $\Omega' =
\Omega$.
\end{proposition}

\proof
By Lemma~\ref{lem:sl2-wheel}, we find
\[
\Omega' = \exp(\sum_n a_{2n} \omega_{2n})
  \equiv \exp(\sum_n 2a_{2n}(\Strutn)^n).
\]
Set $f(x) = \exp(2\sum a_{2n} x^n) = \sum f_n x^n$.  Then by
Lemma~\ref{lem:sl2-gluing},
\begin{align*}
\langle\Omega',(\Strutn)^n\rangle &\equiv
     \langle f(\Strutn), (\Strutn)^n\rangle
    = \langle f_n(\Strutn)^n, (\Strutn)^n\rangle
    \equiv f_n (2n+1)! \\
  &= \left(\frac{1}{24}\ThetaGraph\right)^n \equiv \frac{1}{4^n}.
\end{align*}
so\vspace{-20pt}
\begin{align*}
f_n &= \frac{1}{4^n (2n+1)!} \\
f(x) &= \frac{\sinh(\sqrt{x}/2)}{\sqrt{x}/2} \\
\exp\left(2\sum_n a_n x^n\right) &= \frac{\sinh (x/2)}{x/2} \\
\sum_n a_n x^n &= \frac{1}{2}\log\frac{\sinh(x/2)}{x/2}.\tag*{\qed}
\end{align*}

\end{document}

%% file: gtoutput.tex
\def\ifplaintex{\expandafter\ifx\csname documentclass\endcsname\relax}

\ifplaintex 
\else
\fi

\expandafter\ifx\csname beginpicture\endcsname\relax
\expandafter\ifx\csname documentclass\endcsname\relax
\input pictex \else
\input prepictex \input pictex \input postpictex \fi\fi

\def\gt{{\mathsurround=0pt\it $\cal G\mskip-2mu$eometry \&\ 
$\cal T\!\!$opology}}        

\def\gtp{{\mathsurround=0pt\it $\cal G\mskip-2mu$eometry \&\ 
$\cal T\!\!$opology $\cal P\!$ublications}}  

\def\lognumber#1{\def\thelognumber{#1}}
\def\volumenumber#1{\def\thevolumenumber{#1}}
\def\papernumber#1{\def\thepapernumber{#1}}
\def\volumeyear#1{\def\thevolumeyear{#1}}

\def\pagenumbers#1#2{\def\startpage{#1}\def\finishpage{#2}}
\def\published#1{\def\publishdate{#1}}
\def\proposed#1{\def\theproposer{#1}}
\def\seconded#1{\def\theseconders{#1}}
\def\received#1{\def\receiveddate{#1}}

\def\accepted#1{\def\accepteddate{#1}}

\def\coverauthors#1{\def\thecoverauthors{#1}}
\def\asciiauthors#1{\def\theasciiauthors{#1}}
\def\asciiaddress#1{\def\theasciiaddress{#1}}
\def\asciiemail#1{\def\theasciiemail{#1}}
\long\def\asciiabstract#1{\long\def\theasciiabstract{#1}}

\def\shortauthors#1{\def\theshortauthors{#1}}

\let\\\par\let\thelognumber\relax
\let\thevolumenumber\relax\let\thepapernumber\relax
\let\thevolumeyear\relax\let\thesamplenumber\relax\let\startpage\relax
\let\finishpage\relax\let\publishdate\relax\let\receiveddate\relax
\let\reviseddate\relax\let\accepteddate\relax\let\theasciititle\relax
\let\theasciiauthors\relax\let\theasciiaddress\relax
\let\theasciiabstract\relax
\let\theasciiemail\relax\let\theshortauthors\relax\let\theshorttitle\relax
\let\thecoverauthors\relax

\long\def\maketitlep{   

\count0=\startpage

\gt\hfill      
\beginpicture
\setcoordinatesystem units <0.33truein, 0.33truein> point at 2.2 0.9
\setplotsymbol ({$\cal G$})
\plotsymbolspacing=9truept
\circulararc 315 degrees from 0 1 center at 0 0
\setplotsymbol ({$\cal T$})
\circulararc 315 degrees from 1 -1 center at 1 0
\endpicture
\break
{\small\ifx\thesamplenumber\relax 
Volume \else Sample
\fi\thevolumenumber\ (\thevolumeyear)
\startpage--\finishpage\nl
Published: \publishdate}
\vglue 0.5truein plus 0.4fil minus 0.1truein

{\parskip=0pt\leftskip 0pt plus 1fil\def\\{\par\smallskip}{\ifplaintex\large
\else\Large\fi\bf\thetitle}\par\medskip}   

\vglue 0pt plus 0.1fil 

{\parskip=0pt\leftskip 0pt plus 1fil\def\\{\par}{\sc\theauthors}
\par\medskip}

\vglue 0pt plus 0.1fil 

{\small\parskip=0pt\let\newline\\
{\leftskip 0pt plus 1fil\def\\{\par}{\sl\theaddress}\par}
\expandafter\ifx\theemail\relax    
\relax\else\vglue 5pt plus 0.02fil minus 2pt\def\\{\stdspace{\rm 
and}\stdspace} 
\cl{Email:\stdspace\tt\theemail}\fi
\ifx\theurl\relax                  
\relax\else\vglue 5pt plus 0.02fil minus 2pt\def\\{\stdspace{\rm 
and}\stdspace}
\cl{URL:\stdspace\tt\theurl}\fi\par}

\vglue 7pt plus 0.3fil minus 3pt

{\bf Abstract}
\vglue 5pt plus 0.1fil minus 2pt

\theabstract

\vglue 7pt plus 0.3fil minus 3pt

{\bf AMS Classification numbers}\quad Primary:\quad \theprimaryclass

Secondary:\quad \thesecondaryclass

\vglue 5pt plus 0.3fil minus 2pt

{\bf Keywords}\quad \thekeywords

\vglue 10pt plus 0.5fil minus 5pt

{\small  Proposed: \theproposer\hfill Received: \receiveddate\nl
Seconded: \theseconders\hfill 
\ifx\reviseddate\relax                         
Accepted: \accepteddate                        
\else
Revised: \reviseddate                          
\fi}
\eject
}       

\let\maketitlepage\maketitlep
\let\maketitle\maketitlepage

\font\phead=cmsl9 scaled 950
\font=cmsl9 scaled 1050
\font\pnum=cmbx10 scaled 913
\font\lnum=cmbx10 
\font\pfoot=cmsl9 scaled 950
\font=cmsl9 scaled 1050
\ifplaintex
\headline{\vbox to 0pt{\vskip -4.5mm\line{\small\phead\ifnum
\count0=\startpage ISSN 1364-0380 (on line)
1465-3060 (printed) \hfill {\pnum\folio}\else\ifodd\count0\def\\{ }%
\ifx\theshorttitle\relax\thetitle\else\theshorttitle\fi\hfill{\pnum\folio}
\else\def\\{ and }{\pnum\folio}\hfill\ifx\theshortauthors\relax\theauthors
\else\theshortauthors\fi\fi\fi}\vss}}
\footline{\vbox to 0pt{\vglue 0mm\line{\small\pfoot\ifnum\count0=\startpage
\copyright\ \gtp\hfill\else
\gt, Volume \thevolumenumber\ (\thevolumeyear)\hfill\fi}\vss
}}
\else
\makeatletter
\def\@oddhead{{\small\lhead\ifnum\count0=\startpage ISSN 1364-0380 (on line)
1465-3060 (printed) \hfill {\lnum\number\count0}\else\ifodd\count0
\def\\{ }\ifx\theshorttitle\relax \thetitle \else\theshorttitle\fi\hfill
{\lnum\number\count0}\else\def\\{ and }{\lnum\number\count0}
\hfill\ifx\theshortauthors\relax 
\theauthors\else\theshortauthors\fi\fi\fi}}\def\@evenhead{\@oddhead}
\def\@oddfoot{\small\lfoot\ifnum\count0=\startpage\copyright\ \gtp\hfill\else
\gt, Volume \thevolumenumber\ (\thevolumeyear)\hfill\fi}
\def\@evenfoot{\@oddfoot}
\makeatother
\fi

\newwrite\gtoutfile
\long\gdef\makeheadfile{  
{\def\\{, }\def\s{ }
\immediate\openout\gtoutfile head.xxx
\immediate\write\gtoutfile{To: math@arxiv.org}
\immediate\write\gtoutfile{Subject: put or rep NNNNN:pppp}
\immediate\write\gtoutfile{--text follows this line--}
\immediate\write\gtoutfile{Proxy-for: \ifx\theasciiauthors\relax
\theauthors\else\theasciiauthors\fi\s<\ifx\theasciiemail\relax\theemail\else\theasciiemail\fi>}
\immediate\write\gtoutfile{\noexpand\\}
\immediate\write\gtoutfile{Authors: \ifx\theasciiauthors\relax
\theauthors\else\theasciiauthors\fi}
{\def\\{ }\immediate\write\gtoutfile{Title: \ifx\theasciititle\relax
\thetitle\else\theasciititle\fi}}
\immediate\write\gtoutfile{Subj-class: GT or SG or MG etc}
\immediate\write\gtoutfile{MSC-class: \theprimaryclass\ifx\thesecondaryclass\relax\else, \thesecondaryclass\fi}
\immediate\write\gtoutfile{Journal-ref: Geom. Topol. \thevolumenumber
(\thevolumeyear) \startpage-\finishpage}
\immediate\write\gtoutfile{Comments: Published by Geometry and Topology at}
\immediate\write\gtoutfile{\s\s http://www.maths.warwick.ac.uk/gt/GTVol\thevolumenumber/paper\thepapernumber.abs.html}
\immediate\write\gtoutfile{\noexpand\\}
\immediate\write\gtoutfile{}
\ifx\theasciiabstract\relax
\immediate\write\gtoutfile{\theabstract}\else
\immediate\write\gtoutfile{\theasciiabstract}\fi
\immediate\write\gtoutfile{}
\immediate\write\gtoutfile{\noexpand\\}
\immediate\write\gtoutfile{}
\immediate\closeout\gtoutfile}}  

\def\maketitlepage{\maketitlep\makeheadfile}
\let\maketitle\maketitlepage

%% file: draws/One.tex
\begin{picture}(0,0)%
\includegraphics{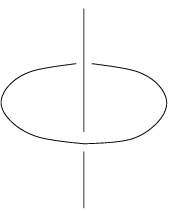}%
\end{picture}%
\setlength{\unitlength}{1657sp}%
\begingroup\makeatletter\ifx\SetFigFont\undefined%
\gdef\SetFigFont#1#2#3#4#5{%
  \reset@font\fontsize{#1}{#2pt}%
  \fontfamily{#3}\fontseries{#4}\fontshape{#5}%
  \selectfont}%
\fi\endgroup%
\begin{picture}(1916,2304)(1113,-2473)
\end{picture}

%% file: draws/OOneTwo.tex
\begin{picture}(0,0)%
\includegraphics{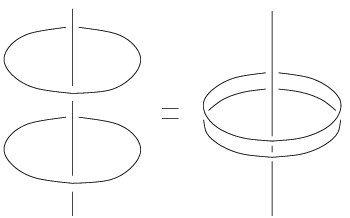}%
\end{picture}%
\setlength{\unitlength}{1367sp}%
\begingroup\makeatletter\ifx\SetFigFont\undefined%
\gdef\SetFigFont#1#2#3#4#5{%
  \reset@font\fontsize{#1}{#2pt}%
  \fontfamily{#3}\fontseries{#4}\fontshape{#5}%
  \selectfont}%
\fi\endgroup%
\begin{picture}(4692,2889)(2883,-3808)
\end{picture}

%% file: draws/HopfSum.tex
\begin{picture}(0,0)%
\includegraphics{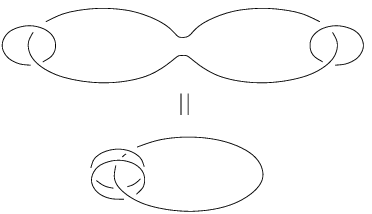}%
\end{picture}%
\setlength{\unitlength}{1492sp}%
\begingroup\makeatletter\ifx\SetFigFont\undefined%
\gdef\SetFigFont#1#2#3#4#5{%
  \reset@font\fontsize{#1}{#2pt}%
  \fontfamily{#3}\fontseries{#4}\fontshape{#5}%
  \selectfont}%
\fi\endgroup%
\begin{picture}(4610,2601)(2354,-3548)
\end{picture}

%% file: draws/Zero.tex
\begin{picture}(0,0)%
\includegraphics{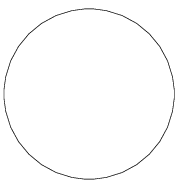}%
\end{picture}%
\setlength{\unitlength}{1036sp}%
\begingroup\makeatletter\ifx\SetFigFont\undefined%
\gdef\SetFigFont#1#2#3#4#5{%
  \reset@font\fontsize{#1}{#2pt}%
  \fontfamily{#3}\fontseries{#4}\fontshape{#5}%
  \selectfont}%
\fi\endgroup%
\begin{picture}(3130,3130)(8129,-5485)
\end{picture}

%% file: draws/nZero.tex
\begin{picture}(0,0)%
\includegraphics{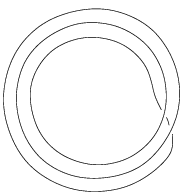}%
\end{picture}%
\setlength{\unitlength}{1036sp}%
\begingroup\makeatletter\ifx\SetFigFont\undefined%
\gdef\SetFigFont#1#2#3#4#5{%
  \reset@font\fontsize{#1}{#2pt}%
  \fontfamily{#3}\fontseries{#4}\fontshape{#5}%
  \selectfont}%
\fi\endgroup%
\begin{picture}(3245,3364)(2275,-4016)
\end{picture}

%% file: draws/A-examp.tex
\begin{picture}(0,0)%
\includegraphics{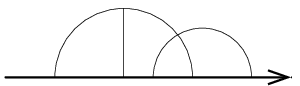}%
\end{picture}%
\setlength{\unitlength}{1243sp}%
\begingroup\makeatletter\ifx\SetFigFont\undefined%
\gdef\SetFigFont#1#2#3#4#5{%
  \reset@font\fontsize{#1}{#2pt}%
  \fontfamily{#3}\fontseries{#4}\fontshape{#5}%
  \selectfont}%
\fi\endgroup%
\begin{picture}(4428,1217)(2699,-1727)
\end{picture}

%% file: draws/B-examp.tex
\begin{picture}(0,0)%
\includegraphics{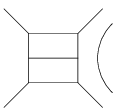}%
\end{picture}%
\setlength{\unitlength}{1036sp}%
\begingroup\makeatletter\ifx\SetFigFont\undefined%
\gdef\SetFigFont#1#2#3#4#5{%
  \reset@font\fontsize{#1}{#2pt}%
  \fontfamily{#3}\fontseries{#4}\fontshape{#5}%
  \selectfont}%
\fi\endgroup%
\begin{picture}(2000,1824)(889,-1873)
\end{picture}

%% file: draws/2wheel.tex
\begin{picture}(0,0)%
\includegraphics{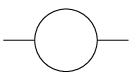}%
\end{picture}%
\setlength{\unitlength}{1973sp}%
\begingroup\makeatletter\ifx\SetFigFont\undefined%
\gdef\SetFigFont#1#2#3#4#5{%
  \reset@font\fontsize{#1}{#2pt}%
  \fontfamily{#3}\fontseries{#4}\fontshape{#5}%
  \selectfont}%
\fi\endgroup%
\begin{picture}(1224,614)(889,-368)
\end{picture}

%% file: draws/4wheel.tex
\begin{picture}(0,0)%
\includegraphics{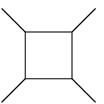}%
\end{picture}%
\setlength{\unitlength}{1973sp}%
\begingroup\makeatletter\ifx\SetFigFont\undefined%
\gdef\SetFigFont#1#2#3#4#5{%
  \reset@font\fontsize{#1}{#2pt}%
  \fontfamily{#3}\fontseries{#4}\fontshape{#5}%
  \selectfont}%
\fi\endgroup%
\begin{picture}(924,924)(3739,-523)
\end{picture}

%% file: draws/6wheel.tex
\begin{picture}(0,0)%
\includegraphics{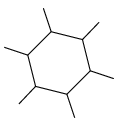}%
\end{picture}%
\setlength{\unitlength}{1973sp}%
\begingroup\makeatletter\ifx\SetFigFont\undefined%
\gdef\SetFigFont#1#2#3#4#5{%
  \reset@font\fontsize{#1}{#2pt}%
  \fontfamily{#3}\fontseries{#4}\fontshape{#5}%
  \selectfont}%
\fi\endgroup%
\begin{picture}(1074,1074)(1864,-1198)
\end{picture}

%% file: draws/SideGlu.tex
\begin{picture}(0,0)%
\includegraphics{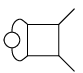}%
\end{picture}%
\setlength{\unitlength}{1973sp}%
\begingroup\makeatletter\ifx\SetFigFont\undefined%
\gdef\SetFigFont#1#2#3#4#5{%
  \reset@font\fontsize{#1}{#2pt}%
  \fontfamily{#3}\fontseries{#4}\fontshape{#5}%
  \selectfont}%
\fi\endgroup%
\begin{picture}(695,624)(366,-223)
\end{picture}

%% file: draws/DiagGlu.tex
\begin{picture}(0,0)%
\includegraphics{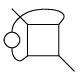}%
\end{picture}%
\setlength{\unitlength}{1973sp}%
\begingroup\makeatletter\ifx\SetFigFont\undefined%
\gdef\SetFigFont#1#2#3#4#5{%
  \reset@font\fontsize{#1}{#2pt}%
  \fontfamily{#3}\fontseries{#4}\fontshape{#5}%
  \selectfont}%
\fi\endgroup%
\begin{picture}(695,624)(666,-523)
\end{picture}